\documentclass{conm-p-l}%
\usepackage{amssymb}
\usepackage{amsmath}
\usepackage{amsfonts}
\usepackage{graphicx}
\newtheorem{theorem}{Theorem}

\newtheorem{conjecture}[theorem]{Conjecture}

\newtheorem{metatheorem}[theorem]{Metatheorem}
\begin{document}

\title{The range of the heat operator}
\author{Brian C. Hall}
\address{University of Notre Dame, Department of Mathematics, 255 Hurley Building, Notre Dame IN, 46556, U.S.A.}
\email{bhall@nd.edu}
\thanks{Supported in part by NSF Grant DMS-0200649}
\date{May 2005}

\begin{abstract}
This paper describes results characterizing the range of the time-$t$ heat
operator on various manifolds, including Euclidean spaces, spheres, and
hyperbolic spaces. The guiding principle behind these results is this: The
functions in the range of the heat operator should be, roughly, those
functions having an analytic continuation to the appropriate complexified
manifold with growth in the imaginary directions at most like that of the
time-$t$ heat heat kernel.

\end{abstract}

\maketitle

\section{Introduction}

The heat equation is very smoothing. Suppose, for example, that we consider
some initial function $f$ in $L^{2}(\mathbb{R}^{d})$ and run the heat equation
for some fixed time $t>0.$ Then the resulting function $F(x):=u(x,t) $
(suppressing the dependence on $t,$ which is fixed) will certainly be
$C^{\infty},$ even real-analytic. Is there some regularity condition that
characterizes which functions $F$ arise in this way? If we define the time-$t
$ heat operator to be the operator taking $f$ to $F,$ we may formulate the
question in this way: Which functions are in the range of the time-$t$ heat
operator? To put it a different way, for which initial conditions $F$ can we
solve the \textit{backward} heat equation for time $t$? Clearly, $F$ must be
real-analytic, but this is not nearly sufficient. The question also makes
sense for the heat equation on other manifolds, with a similar situation:
Functions in the range of the heat operator are real-analytic, but this is not
sufficient to characterize the image.

One approach to characterizing the range of the heat operator is to use the
Fourier transform. Let us consider initial conditions in $L^{2}(\mathbb{R}%
^{d})$ and let us normalize the heat equation as $\partial u/\partial
t=\frac{1}{2}\Delta u.$ Then for any one fixed $t>0,$ the range of the
time-$t$ heat operator consists of precisely those functions $F$ such that%
\[
\int_{\mathbb{R}^{d}}\left\vert \hat{F}(\xi)\right\vert ^{2}e^{t\left\vert
\xi\right\vert ^{2}}d\xi<\infty,
\]
where $\hat{F}$ denotes the Fourier transform of $F.$ Although this is a very
simple condition, it is desirable to have a condition expressible in terms of
$F$ itself, rather than in terms of the Fourier transform of $F.$ On other
manifolds, similarly, one could characterize the image of the heat operator in
terms of the expansion of a function in eigenfunctions of the Laplacian.
Nevertheless, it would be desirable to have a characterization in terms of the
function itself.

The premise of this paper is that the right way to characterize the range of
the heat operator (in terms of the functions themselves) is by means of
analyticity properties. That is, I wish to formulate the question as: How
analytic are the functions in the range of the heat operator? The answer to
this question should be provided by the (ubiquitous) heat kernel. The heat
operator is computed by integration against the heat kernel, so functions in
the range of the time-$t$ heat operator should pick up the analyticity
properties of the time-$t$ heat kernel. This paper, then, will investigate the
validity of the following (deliberately vague) idea.

\begin{quote}
\noindent For any one fixed positive time $t,$ the functions in the range of
the time-$t$ heat operator are the functions having the same analyticity
properties as the time-$t$ heat kernel.
\end{quote}

I do not know how to make this idea work for general manifolds. I will
consider mainly certain special manifolds, the Euclidean spaces, spheres,
hyperbolic spaces and other compact and noncompact symmetric spaces. The cases
of Euclidean space and compact symmetric spaces (including spheres) are by now
well understood. The case of a noncompact symmetric space is just beginning to
be understood, and I\ concentrate on the most tractable example, hyperbolic 3-space.

Let us now be a bit more precise about how we will try to use analyticity
properties to characterize the range of the heat operator. For each manifold
$\mathcal{M}$ that we consider, we will introduce an appropriate
\textquotedblleft complexification\textquotedblright\ $\mathcal{M}%
_{\mathbb{C}}.$ The time-$t$ heat kernel for $\mathcal{M}$ will then have an
analytic continuation to $\mathcal{M}_{\mathbb{C}}$ with a certain
$t$-dependent rate of growth in the \textquotedblleft imaginary
directions.\textquotedblright\ The rate of growth of the analytically
continued heat kernel will lead to a \textquotedblleft
metatheorem\textquotedblright\ characterizing, roughly, the range of the
time-$t$ heat operator as those functions having analytic continuations to
$\mathcal{M}_{\mathbb{C}}$ with at most the specified growth rate in the
imaginary direction. For each manifold, the metatheorem will be translated
into various precisely stated theorems characterizing the image under the
time-$t$ heat operator of various spaces of initial conditions.

The strategy in the previous paragraph can be implemented in a straightforward
way for Euclidean spaces, spheres, and compact symmetric spaces. The case of
hyperbolic space or a noncompact symmetric space is substantially more subtle,
because of singularities that arise in the analytically continued heat kernel.
Nevertheless, recent results allow for progress in these cases, as described
in Section \ref{noncompact.sec}. Section 6 gives a brief account of the issues
that are involved in the case of general manifolds.

I thank the referee for making useful corrections to the manuscript.

\section{The $\mathbb{R}^{d}$ case}

We consider the Laplacian $\Delta$ on $\mathbb{R}^{d},$ given by%
\[
\Delta=\sum_{k=1}^{d}\frac{\partial^{2}}{\partial x_{k}^{2}}.
\]
Note that our Laplacian is a \textit{negative} operator, which means that
$e^{t\Delta/2}$ is the \textit{forward} heat operator. We may view $\Delta$ as
an unbounded operator in $L^{2}(\mathbb{R}^{d},dx),$ with domain initially
taken to be, say, $C^{\infty}$ functions of compact support. The Laplacian is
essentially self-adjoint on this domain. We let $\Delta$ also denote the
closure of this operator, which is self-adjoint and negative-semidefinite. The
spectral theorem then allows us to define the heat operator $e^{t\Delta/2}$ as
a contraction operator on $L^{2}(\mathbb{R}^{d}).$ We will let
\textquotedblleft time-$t$ heat operator\textquotedblright\ denote
$e^{t\Delta/2}$ (with a factor of 2 in the exponent).

Given $f$ in $L^{2}(\mathbb{R}^{d},dx),$ if we define%
\[
u(x,s)=(e^{s\Delta/2}f)(x),
\]
then $u$ satisfies the heat equation%
\[
\frac{\partial u}{\partial s}=\frac{1}{2}\Delta u,\quad s>0,
\]
subject to the initial condition%
\[
\lim_{s\searrow0}u(x,s)=f(x).
\]
Here the limit is in the norm topology of $L^{2}(\mathbb{R}^{d},dx)$.

It is well known that the heat operator can be computed by convolution against
the \textit{heat kernel}, namely, the function%
\begin{equation}
\rho_{t}(x)=(2\pi t)^{-d/2}e^{-x^{2}/2t},\label{rho.def1}%
\end{equation}
where $x^{2}=x_{1}^{2}+\cdots+x_{d}^{2}.$ This means that%
\begin{equation}
(e^{t\Delta/2}f)(x)=\int_{\mathbb{R}^{d}}\rho_{t}(x-x^{\prime})f(x^{\prime
})~dx^{\prime}.\label{convolution}%
\end{equation}
In this expression, we initially think of $x$ as being in $\mathbb{R}^{d}.$
However, it is easy to see that the heat kernel (\ref{rho.def1}) has an entire
analytic continuation to $\mathbb{C}^{d}.$ Furthermore, on the right-hand side
of (\ref{convolution}), the variable $x$ appears only in $\rho_{t}$ and not in
$f.$ Thus, even if $f$ itself is not analytic, we can analytically continue
$e^{t\Delta/2}f$ from $\mathbb{R}^{d}$ to $\mathbb{C}^{d}$ by setting%
\begin{align}
(e^{t\Delta/2}f)(z)  & =\int_{\mathbb{R}^{d}}\rho_{t}(z-x^{\prime}%
)f(x^{\prime})~dx^{\prime}\nonumber\\
& =(2\pi t)^{-d/2}\int_{\mathbb{R}^{d}}e^{-(z-x^{\prime})^{2}/2t}f(x^{\prime
})~dx^{\prime},\quad z\in\mathbb{C}^{d}.\label{convolution2}%
\end{align}
Here $(z-x^{\prime})^{2}$ is to be interpreted as $(z_{1}-x_{1}^{\prime}%
)^{2}+\cdots+(z_{d}-x_{d}^{\prime})^{2}$ (no absolute values). It is not hard
to prove, using Morera's Theorem, that this integral really depends
holomorphically on $z,$ for any $f$ in $L^{2}(\mathbb{R}^{d}).$

Note that the analytically continued heat kernel $\rho_{t}(z)$ grows like
$e^{\left\vert y\right\vert ^{2}/2t}$ in the imaginary directions.
Specifically, if $z=x+iy,$ with $x$ and $y$ in $\mathbb{R}^{d},$ then%
\begin{equation}
\left\vert \rho_{t}(z-x^{\prime})\right\vert =(2\pi t)^{-d/2}e^{-|x-x^{\prime
}|^{2}/2t}e^{\left\vert y\right\vert ^{2}/2t}.\label{rhot.abs}%
\end{equation}
Suppose we take a function $f,$ convolve $f$ with $\rho_{t},$ and then
analytically continue $e^{t\Delta/2}f$ to $\mathbb{C}^{d},$ as in
(\ref{convolution2}). Then the analytic continuation of $e^{t\Delta/2}f$ will
have at most the same rate of growth in the imaginary directions as $\rho_{t}$
itself. This suggests the following \textquotedblleft
metatheorem,\textquotedblright\ which will serve as the unifying idea behind
the theorems stated in this section. As usual, we are considering the heat
operator for one \textit{fixed} positive time $t.$

\begin{metatheorem}
The functions in the range of the heat operator $e^{t\Delta/2}$ for
$\mathbb{R}^{d}$ are those functions $F$ having an analytic continuation to
$\mathbb{C}^{d}$ with growth at most like $e^{\left\vert y\right\vert ^{2}/2t}
$ in the imaginary directions.
\end{metatheorem}

This metatheorem is worded in a deliberately vague way and is not intended to
be taken too literally. In the first place, the metatheorem does not say what
sort of functions we take in the \textit{domain} of the heat operator. We will
consider various spaces of initial conditions: $L^{2}(\mathbb{R}^{d}) $ with
respect to Lebesgue measure, $L^{2}(\mathbb{R}^{d})$ with respect to various
Gaussian measures, Sobolev spaces, the Schwarz space, and the space of
tempered distributions. In the second place, the metatheorem does not say
precisely how \textquotedblleft growth at most like $e^{\left\vert
y\right\vert ^{2}/t}$ in the imaginary directions\textquotedblright\ is to be
interpreted. We will be flexible in how we interpret this condition, sometimes
in an $L^{2}$ sense, sometimes in a pointwise sense, and in either case with
different possible dependence of the estimates on the real variable $x.$ We
will also polynomial variations around the \textquotedblleft
basic\textquotedblright\ growth rate of $e^{\left\vert y\right\vert ^{2}/t}$
to account for different levels of smoothness in the initial data. I hope
that, despite all this equivocation, each of the precisely formulated theorems
in this section can be seen to be compatible with our metatheorem.

\subsection{Isometry formulas}

In this subsection, we consider mainly initial data in $L^{2}(\mathbb{R}%
^{d}),$ either with respect to Lebesgue measure or with respect to a Gaussian
measure. This case is the easiest to deal with because it allows the use of
Hilbert space methods. This case is also natural from the point of view of
quantum mechanics, where the heat equation arises in the so-called
Segal--Bargmann transform \cite{Se1,Se2,Se3,Ba1}.

We now come to our first theorem, characterizing the image of $L^{2}%
(\mathbb{R}^{d},dx)$ under the time-$t$ heat operator, where $dx$ denotes
Lebesgue measure on $\mathbb{R}^{d}.$ This result is a slight modification of
results of Segal \cite{Se3} and Bargmann \cite{Ba1}. See \cite[Sect.
6.3]{mexnotes} or \cite{newform} for how to adapt the results of Segal and
Bargmann to the present setting.

\begin{theorem}
[Isometry 1]\label{isometry.thm1}A function $F$ on $\mathbb{R}^{d}$ is of the
form $e^{t\Delta/2}f,$ with $f$ in $L^{2}(\mathbb{R}^{d},dx),$ if and only if
$F$ has an analytic continuation to $\mathbb{C}^{d}$ with the property that%
\begin{equation}
\int_{\mathbb{C}^{d}}\left\vert F(x+iy)\right\vert ^{2}\frac{e^{-\left\vert
y\right\vert ^{2}/t}}{(\pi t)^{d/2}}~dy~dx<\infty.\label{isometry.int1a}%
\end{equation}
If $F$ is such a function then%
\begin{equation}
\int_{\mathbb{C}^{d}}\left\vert F(x+iy)\right\vert ^{2}\frac{e^{-\left\vert
y\right\vert ^{2}/t}}{(\pi t)^{d/2}}~dy~dx=\int_{\mathbb{R}^{d}}\left\vert
f(x)\right\vert ^{2}dx.\label{isometry.int1b}%
\end{equation}

\end{theorem}

The isometry theorem may be restated as saying that $F$ is in the range of the
time-$t$ heat operator if and only if $F(x+iy)e^{-\left\vert y\right\vert
^{2}/2t}$ is square-integrable with respect to Lebesgue measure on
$\mathbb{C}^{d}.$ This condition is just an $L^{2}$ version of the growth
condition in our metatheorem. That is to say, in Theorem \ref{isometry.thm1},
\textquotedblleft$F$ has growth at most like $e^{\left\vert y\right\vert
^{2}/2t}$\textquotedblright\ is interpreted to mean that $F(x+iy)$ divided by
$e^{\left\vert y\right\vert ^{2}/2t}$ is square-integrable with respect to
Lebesgue measure.

It is instructive to see how this $L^{2}$ growth condition relates to
pointwise growth conditions. By adapting a result of Bargmann \cite{Ba1} to
our normalization of the isometry theorem (compare \cite[Thm. 6]{H1}) we see
that if a holomorphic function $F$ on $\mathbb{C}^{d}$ satisfies
(\ref{isometry.int1a}), then $F$ automatically satisfies the pointwise bound%
\begin{equation}
\left\vert F(x+iy)\right\vert \leq Ce^{\left\vert y\right\vert ^{2}%
/2t}\label{bound1}%
\end{equation}
where%
\[
C=\left\Vert f\right\Vert _{L^{2}(\mathbb{R}^{d},dx)}(4\pi t)^{d/4}.
\]
Meanwhile, direct calculation shows that if $F$ satisfies a polynomially
stronger bound, say,%
\begin{equation}
\left\vert F(x+iy)\right\vert \leq\frac{Ce^{\left\vert y\right\vert ^{2}/2t}%
}{(1+\left\vert x\right\vert ^{2}+\left\vert y\right\vert ^{2}%
)^{d/2+\varepsilon}},\label{suff.bound1}%
\end{equation}
then $F$ satisfies (\ref{isometry.int1a}). Thus, the $L^{2}$ version of
\textquotedblleft growth at most like $e^{\left\vert y\right\vert ^{2}/2t}%
$\textquotedblright\ in (\ref{isometry.int1a}) is closely related to (but not
equivalent to) the pointwise version of this condition in (\ref{bound1}).

We could also consider initial conditions that are square integrable with
respect to some measure besides Lebesgue measure. The most important such
measures from the point of view of physics (and also the infinite-dimensional
case) are the Gaussian measures. We consider, then, $L^{2}(\mathbb{R}^{d}%
,\rho_{s}(x)dx),$ where $\rho_{s}$ is the Gaussian density defined in
(\ref{rho.def1}). A simple calculation shows that the convolution
(\ref{convolution}) is convergent for all $f$ in $L^{2}(\mathbb{R}^{d}%
,\rho_{s}(x)dx),$ provided that $t<2s.$ We now state a result from
\cite{DH1,newform}.

\begin{theorem}
[Isometry 2]\label{isometry.thm2}For $t<2s,$ a function $F$ on $\mathbb{R}%
^{d}$ is of the form $F=e^{t\Delta/2}f$, with $f\in L^{2}(\mathbb{R}^{d}%
,\rho_{s}(x)dx),$ if and only if $F$ has an analytic continuation to
$\mathbb{C}^{d}$ satisfying%
\begin{equation}
\int_{\mathbb{C}^{d}}\left\vert F(x+iy)\right\vert ^{2}\frac{e^{-\left\vert
y\right\vert ^{2}/t}}{(\pi t)^{d/2}}\frac{e^{-\left\vert x\right\vert ^{2}/r}%
}{(\pi r)^{d/2}}~dy~dx<\infty,\label{isometry.int2a}%
\end{equation}
where%
\[
r=2s-t.
\]
If $F$ is such a function then%
\begin{equation}
\int_{\mathbb{C}^{d}}\left\vert F(x+iy)\right\vert ^{2}\frac{e^{-\left\vert
y\right\vert ^{2}/t}}{(\pi t)^{d/2}}\frac{e^{-\left\vert x\right\vert ^{2}/r}%
}{(\pi r)^{d/2}}~dy~dx=\int_{\mathbb{R}^{d}}\left\vert f(x)\right\vert
^{2}\frac{e^{-\left\vert x\right\vert ^{2}/2s}}{(2\pi s)^{d/2}}%
~dx.\label{isometry.int2b}%
\end{equation}

\end{theorem}

The condition (\ref{isometry.int2a}) is another $L^{2}$ version of the
condition \textquotedblleft growth like $e^{\left\vert y\right\vert ^{2}/2t}$
in the imaginary directions,\textquotedblright\ except that we now allow the
constant in that growth rate to depend on $x.$ The growth rate in the
$x$-direction, namely, $e^{\left\vert x\right\vert ^{2}/2r},$ is what one
would expect based on estimating the rate of growth of $(e^{t\Delta/2}f)(x)$
for $x\in\mathbb{R}^{d}$ and $f\in L^{2}(\mathbb{R}^{d},\rho_{s}(x)dx).$
(Compare Section 2 of \cite{lpbounds}.)

Using the methods of \cite[Sect. 9]{H1}, one can compute the reproducing
kernel for the holomorphic $L^{2}$ space in (\ref{isometry.int2a}) (i.e., the
space of holomorphic functions satisfying (\ref{isometry.int2a})). The result
is that if $F$ is holomorphic and satisfies (\ref{isometry.int2a}) then $F$
satisfies a pointwise bound of the form%
\begin{equation}
\left\vert F(x+iy)\right\vert \leq Ce^{\left\vert y\right\vert ^{2}%
/2t}e^{\left\vert x\right\vert ^{2}/2r}.\label{bound2}%
\end{equation}
Conversely, a holomorphic function satisfying a polynomially stronger bound
than this will, by direct calculation, satisfy (\ref{isometry.int2a}). Thus,
as in the case of the image of $L^{2}(\mathbb{R}^{d},dx),$ the $L^{2}$
condition in (\ref{isometry.int2a}) is closely related to, but not identical
to, the analogous pointwise condition.

The case $s=t$ is just a standard form of the Segal--Bargmann transform
(essentially just the finite-dimensional version of the transform in
\cite{Se3} or \cite{BSZ}). The general case is developed in \cite{DH1} and
\cite{newform}. The case of $L^{2}(\mathbb{R}^{d},dx)$ can be obtained, at
least heuristically, by rescaling the norms on both sides of
(\ref{isometry.int2b}) by a factor of $(2\pi s)^{d/2}$ and then letting $s$
tend to infinity. In the other direction, if we allow $s$ to approach $t/2,$
then the quantity $r=2s-t$ tends to zero. This means that as $s$ approaches
$t/2,$ the measure on $\mathbb{C}^{d}$ collapses onto the imaginary axis. In
this limit we obtain a finite-dimensional version of the Fourier--Wiener transform.

Finally, let us consider briefly initial data in $L^{p}(\mathbb{R}^{d},dx),$
with $1\leq p\leq\infty.$ (See also \cite{lpbounds}.) For $p\neq2,$ we take
(\ref{convolution}) as the definition of $e^{t\Delta/2}.$ It is unlikely that
there is any simple exact description of the image of the heat operator for
$p\neq2.$ Nevertheless, we can give a necessary condition and a sufficient
condition for $F$ to be in the image of $L^{p}(\mathbb{R}^{d},dx),$ with not
too much of a gap between the two conditions. From (\ref{rhot.abs}) we see
that for each fixed $z\in\mathbb{C}^{d}$ we have%
\[
\left\Vert \rho_{t}(z-x^{\prime})\right\Vert _{L^{p^{\prime}}(\mathbb{R}%
^{d},dx^{\prime})}=Ce^{\left\vert \operatorname{Re}z\right\vert ^{2}/2t},
\]
where $p^{\prime}$ is the conjugate exponent to $p$ and where $C$ is a
constant depending on $p,$ $d,$ and $t,$ but not on $z.$ From this it follows
that if $F=e^{t\Delta/2}f$ with $f\in L^{p}(\mathbb{R}^{d},dx)$ we have%
\begin{equation}
\left\vert F(x+iy)\right\vert \leq C\left\Vert f\right\Vert _{L^{p}%
(\mathbb{R}^{d},dx)}e^{\left\vert y\right\vert ^{2}/2t}.\label{lp.bnd1}%
\end{equation}

Conversely, using the inversion formula (\ref{inverse.adjoint}) from Section
\ref{rdinv.sect} and moving the $L^{p}$ norm inside the integral, we can show
that if%
\begin{equation}
\int_{\mathbb{C}^{d}}\left\vert F(x+iy)\right\vert e^{-\left\vert y\right\vert
^{2}/2t}~dx~dy<\infty,\label{lp.bnd2}%
\end{equation}
then there exists a function $f$ on $\mathbb{R}^{d}$ such that $F=e^{t\Delta
/2}f$ and such that $f\in L^{p}(\mathbb{R}^{d},dx)$ for all $1\leq p\leq
\infty.$ We may say, then, that the image of $L^{p}(\mathbb{R}^{d},dx)$ is
somewhere between \textquotedblleft growth like $e^{\left\vert y\right\vert
^{2}/2t}~$\textquotedblright\ interpreted in the pointwise sense (condition
(\ref{lp.bnd1})) and \textquotedblleft growth like $e^{\left\vert y\right\vert
^{2}/2t}~$\textquotedblright\ interpreted in the $L^{1}$ sense (condition
(\ref{lp.bnd2})).

By a simple interpolation argument (as in \cite{lpbounds}) we can improve the
results in the previous two paragraphs as follows. For $1\leq p\leq2,$ if $f$
is in $L^{p}(\mathbb{R}^{d},dx)$ then $Fe^{-\left\vert y\right\vert ^{2}/2t}$
is in $L^{p^{\prime}}(\mathbb{C}^{d},dz)$. For $2\leq p\leq\infty,$ if
$Fe^{-\left\vert y\right\vert ^{2}/2t}$ is in $L^{p^{\prime}}(\mathbb{C}%
^{d},dz),$ then $f$ is in $L^{p}(\mathbb{R}^{d},dx).$

\subsection{Inversion formula\label{rdinv.sect}}

In this section, we consider the process of recovering $f$ from (the analytic
continuation of) $F:=e^{t\Delta/2}f.$ The inversion formula in Theorem
\ref{inv.thm1} will also give another piece of evidence for the validity of
our metatheorem. Specifically, the inversion formula will help explain why a
growth rate of $e^{\left\vert y\right\vert ^{2}/2t}$ should be sufficient to
imply that $F$ is in the range of the heat operator. (Analysis of the
convolution (\ref{convolution}) shows that this condition is necessary;
sufficiency is not as obvious.) Of course, the isometry theorems of the
previous subsection already tell us that the growth rate of $e^{\left\vert
y\right\vert ^{2}/2t}$ (interpreted in an $L^{2}$ sense) is sufficient.
Nevertheless, the inversion formula will give us a more intuitive reason for
sufficiency: A growth rate of $e^{\left\vert y\right\vert ^{2}/2t}$ allows a
direct construction of the function $f$ for which $F=e^{t\Delta/2}f.$

It is important to note that there are many different possible inversion
formulas for the heat operator. After all, since $F$ is a holomorphic
function, there can be many different integrals involving $F$ that all give
the same answer. For example, in the Cauchy integral formula, integrals over
many different contours all yield the value $F(z).$ See \cite{lpbounds} for a
application of multiple inversion formulas for the heat operator.

The most obvious inversion formula is the one obtained by thinking of the heat
operator (followed by analytic continuation) as a unitary map from
$L^{2}(\mathbb{R}^{d},dx)$ onto the Hilbert space of $L^{2}$ holomorphic
functions on $\mathbb{C}^{d}$ with respect to the Gaussian measure in
(\ref{isometry.int1a}). Since the adjoint of a unitary map is its inverse and
since the transform may be thought of as an integral operator with integral
kernel $\rho_{t}(z-x),$ we easily obtain the following inversion formula,
valid for any $f\in L^{2}(\mathbb{R}^{d},dx)$:%
\begin{equation}
f(x)=\lim_{R\rightarrow\infty}\int_{\left\vert z\right\vert \leq R}%
\overline{\rho_{t}(z-x)}F(z)\frac{e^{-\left\vert \operatorname{Re}z\right\vert
^{2}/t}}{(\pi t)^{d/2}}~dz,\label{inverse.adjoint}%
\end{equation}
where the limit is in the norm topology of $L^{2}(\mathbb{R}^{d},dx).$
(Compare Section 9 of \cite{H1}.)

It is tempting to think of (\ref{inverse.adjoint}) as \textit{the} inversion
formula for the heat operator. However, there is another inversion formula
that is, in some ways, more fundamental.

\begin{theorem}
[Inversion]\label{inv.thm1}Let $f$ be in $L^{2}(\mathbb{R}^{d},dx)$ and let
$F=e^{t\Delta/2}f.$ Then $F$ has an analytic continuation to $\mathbb{C}^{d}$
and $f$ may be recovered from this analytic continuation by the formula%
\begin{equation}
f(x)=\lim_{R\rightarrow\infty}\int_{\left\vert y\right\vert \leq
R}F(x+iy)\frac{e^{-\left\vert y\right\vert ^{2}/2t}}{(2\pi t)^{d/2}%
}~dy,\label{inv.int1}%
\end{equation}
where the limit is in the norm topology of $L^{2}(\mathbb{R}^{d},dx).$
\end{theorem}

This inversion formula is special because the value of $f$ at a point $x$ is
computed from the values of $F$ at points of the form $x+iy$ (at least in
those cases where the limit exists pointwise). This inversion formula is also
special because it is easy (as we shall see) to prove directly that the
integral in (\ref{inv.int1}) undoes the heat equation. The analog of
(\ref{inv.int1}) in the case of compact symmetric spaces (Sections
\ref{sphere.sec} and \ref{compact.sec} below) is an important part of the
theory and helps elucidate the role of duality between compact and noncompact
symmetric spaces.

In general, one cannot expect pointwise convergence of the limit in
(\ref{inv.int1}). After all, $f$ is an arbitrary $L^{2}$ function, which can
have singularities. Any inversion formula must have a mechanism for
reproducing those singularities, and in (\ref{inv.int1}) that mechanism is the
possible failure of convergence of the integral. An alternative means of
regularizing the integral is
\begin{equation}
f(x)=\lim_{s\nearrow t}\int_{\mathbb{R}^{d}}F(x+iy)\frac{e^{-\left\vert
y\right\vert ^{2}/2s}}{(2\pi s)^{d/2}}~dy,\label{inv.int2}%
\end{equation}
where the pointwise bounds (\ref{bound1}) show that the integral
(\ref{inv.int2}) is absolutely convergent for all $s<t.$ Here, again, the
limit is to be taken in the norm topology of $L^{2}(\mathbb{R}^{d},dx).$ As we
will see in the next subsection, if we assume sufficient smoothness for $f,$
then we may take $s=t$ in (\ref{inv.int2}) (or $R=\infty$ in (\ref{inv.int1}%
)), with absolute convergence of the resulting integral.

Let us think about why the inversion formula works. We give now an intuitive
argument; a rigorous proof follows. The key observation is that the function%
\begin{equation}
\frac{e^{-\left\vert y\right\vert ^{2}/2s}}{(2\pi s)^{d/2}}\label{y.gauss}%
\end{equation}
occurring in (\ref{inv.int2}) is just the heat kernel for $\mathbb{R}^{d}$
(compare (\ref{rho.def1})), except now thought of a living on the imaginary
axis in $\mathbb{C}^{d}.$ The integral (\ref{inv.int1}) is then effectively
computing the (forward) heat operator in the $y$-directions. However, because
$F$ is holomorphic, the Laplacian in the $x$-direction is the negative of the
Laplacian in the $y$-direction. Thus, doing the \textit{forward} heat operator
in the $y$-variables amounts to doing the \textit{backward} heat operator in
the $x$-variables, thus undoing the heat equation.

More precisely, let us fix a holomorphic function $F$ and define%
\[
u(x,s)=\int_{\mathbb{R}^{d}}F(x+iy)\frac{e^{-\left\vert y\right\vert ^{2}/2s}%
}{(2\pi s)^{d/2}}~dy
\]
whenever the integral converges. Assume the integral has reasonable
convergence properties at least for $s<t.$ Then the following argument shows
that $u$ satisfies the \textit{backward} heat equation on $\mathbb{R}%
^{d}\times(0,t)$: (1) differentiate under the integral with respect to $s$;
(2) change the $s$-derivative on the Gaussian factor (\ref{y.gauss}) into
$\frac{1}{2}\partial^{2}/\partial y^{2}$; (3) integrate by parts to move the
$y$-derivatives onto $F$; (4) use the Cauchy--Riemann equations for $F$ to
change $\partial^{2}/\partial y^{2}$ into $-\partial^{2}/\partial x^{2},$
which then comes outside the integral. Meanwhile, as $s$ approaches zero,
$u(x,s)$ will converge to $F(x).$ If we assume also that $u(x,s)$ approaches
some function $f(x)$ as $s$ approaches $t,$ then we will have that
$f=e^{-t\Delta/2}F$ and $F=e^{t\Delta/2}f.$ Thus (\ref{inv.int2}) indeed
serves to recover $f$ from $F.$

Now, the above argument also helps us see how the inversion formula fits into
our metatheorem. What conditions should $F$ satisfy so that we can run the
above argument and produce a function $f$ with $F=e^{t\Delta/2}f$? We need
convergence of the integral (\ref{inv.int2}) for $s<t$ and reasonable behavior
of the integral as $s$ approaches $t.$ That is, we need (roughly) $F $ to have
growth at most like $e^{\left\vert y\right\vert ^{2}/2t}$ in the imaginary
directions. Thus, even if this reasoning does not easily translate into a
precise \textquotedblleft if and only if\textquotedblright\ result, the
inversion formula still gives another confirmation of the validity of our metatheorem.

Note that in the isometry formula (\ref{inv.int1}) we have $e^{-\left\vert
y\right\vert ^{2}/t}$, whereas in the inversion formula we have
$e^{-\left\vert y\right\vert ^{2}/2t}.$ This factor of two change in the
variance is necessary because the isometry formula involves $\left\vert
F\right\vert ^{2}$ (which grows like $e^{\left\vert y\right\vert ^{2}/t}$)
whereas the inversion formula involves $F$ itself (which grows like
$e^{\left\vert y\right\vert ^{2}/2t}$).

\begin{proof}[Proof of Theorem \ref{inv.thm1}] The problem with the direct argument above for the
inversion formula is that it does not easily yield convergence in $L^{2}$ as
$R\rightarrow\infty.$ We use, then, another argument involving the Fourier
transform. Let us normalize the Fourier transform in such a way that the
Fourier inversion formula takes the form $f(x)=\int e^{i\xi x}\hat{f}(\xi
)d\xi.$ A simple calculation then shows that%
\begin{align}
& \int_{\left\vert y\right\vert \leq R}F(x+iy)\frac{e^{-\left\vert
y\right\vert ^{2}/2t}}{(2\pi t)^{d/2}}~dy\label{inv.calc}\\
& =\int_{\mathbb{R}^{d}}e^{i\xi x}~\left[  \hat{f}(\xi)\left(  \int
_{\left\vert y\right\vert \leq R}e^{-\xi y}\frac{e^{-\left\vert y\right\vert
^{2}/2t}}{(2\pi t)^{d/2}}~dy\right)  e^{-\left\vert \xi\right\vert ^{2}%
t/2}\right]  ~d\xi.\nonumber
\end{align}
Now, the expression in round parentheses converges monotonically to
$e^{\left\vert \xi\right\vert ^{2}t/2}$ as $R$ tends to infinity. It follows
that the expression in square brackets is in $L^{2}$ for all $R$ and converges
to $\hat{f}(\xi)$ in $L^{2}.$ Thus, the left-hand side of (\ref{inv.calc})
converges in $L^{2}$ to $f$ as $R$ tends to infinity.
\end{proof}

\subsection{The effect of smoothness of the initial data}

Up to now, we have considered the behavior of $F:=e^{t\Delta/2}f$ when $f$ is
in some $L^{2}$ or $L^{p}$ space. It is worth considering how the behavior of
$F$ will change if we assume some differentiability for $f.$ Our main result
will concern the image under $e^{t\Delta/2}$ of the Sobolev space
$H^{n}(\mathbb{R}^{d}),$ which consists of functions whose partial derivatives
up to order $n$ (computed in the distributional sense) lie in $L^{2}%
(\mathbb{R}^{d},dx).$ We will also record results of Bargmann concerning the
image under $e^{t\Delta/2}$ of the Schwarz space and of the space of tempered distributions.

\begin{metatheorem}
Smoothness properties of $f$ are reflected in polynomial fluctuations of $F$
around the \textquotedblleft basic\textquotedblright\ growth rate of
$e^{\left\vert y\right\vert ^{2}/2t}$ in the imaginary directions.
Specifically, each derivative that $f$ has in $L^{2}$ gives an improvement in
the behavior of $F(x+iy)$ by one power of $\left\vert y\right\vert $.
\end{metatheorem}

This gets interpreted, as in the isometry theorems, in an $L^{2}$ sense.

\begin{theorem}
\label{sobolev.thm1}Suppose $n$ is a positive integer. Then a function $F$ on
$\mathbb{R}^{d}$ is of the form $F=e^{t\Delta/2}f,$ with $f$ in the Sobolev
space $H^{n}(\mathbb{R}^{d}),$ if and only if $F$ has an analytic continuation
to $\mathbb{C}^{d}$ satisfying%
\begin{equation}
\int_{\mathbb{C}^{d}}\left\vert F(x+iy)\right\vert ^{2}e^{-\left\vert
y\right\vert ^{2}/t}(1+\left\vert y\right\vert ^{2n})~dy~dx<\infty
.\label{sobolev.int1}%
\end{equation}
If $F$ is such a function then $F$ satisfies the pointwise bounds%
\begin{equation}
\left\vert F(x+iy)\right\vert \leq C\frac{e^{\left\vert y\right\vert ^{2}/2t}%
}{1+\left\vert y\right\vert ^{n}}\label{sobolev.bnd1}%
\end{equation}
for some constant $C.$
\end{theorem}

We are \textit{not} asserting that the norm in $H^{n}$ is equal to the
integral in (\ref{sobolev.int1}). Rather, the norm in (\ref{sobolev.int1}) is
different from but equivalent to the $H^{n}$ norm on $f.$ (Compare
(\ref{sob.norm}) in the proof of Theorem \ref{sobolev.thm1}.)

Note that if $n>d,$ then (\ref{sobolev.bnd1}) implies that $\int
_{\mathbb{R}^{d}}\left\vert F(x+iy)\right\vert e^{-\left\vert y\right\vert
^{2}/2t}~dy<\infty$. This means that we may take $R=\infty$ in (\ref{inv.int1}%
) and obtain%
\begin{equation}
f(x)=\int_{\mathbb{R}^{d}}F(x+iy)\frac{e^{-\left\vert y\right\vert ^{2}/2t}%
}{(2\pi t)^{d/2}}~dy,\label{inv.sob}%
\end{equation}
with absolute convergence of the integral for all $x.$ However, the condition
$n>d$ is not sharp. By a better argument (which does not use pointwise bounds
as an intermediate step) we can get the same result whenever $n>d/2.$ This is
presumably the optimal condition, since this is also the condition in the
Sobolev embedding theorem needed to guarantee that $f\in H^{n}(\mathbb{R}%
^{d})$ is continuous. (If $n\leq d/2$ then $f\in H^{n}(\mathbb{R}^{d})$ can
have singularities, in which case we cannot expect convergence of the integral
(\ref{inv.sob}) for all $x.$)

\begin{theorem}
\label{pointwise.thm}If $f\in H^{n}(\mathbb{R}^{d})$ and $n>d/2,$ then for all
$x\in\mathbb{R}^{d}$ we have%
\[
f(x)=\int_{\mathbb{R}^{d}}F(x+iy)\frac{e^{-\left\vert y\right\vert ^{2}/2t}%
}{(2\pi t)^{d/2}}~dy
\]
with absolute convergence of the integral for all $x.$
\end{theorem}

\begin{proof}[Proof of Theorem \ref{sobolev.thm1}] I adapt here the arguments in \cite{HL} from
the compact group case to the (easier) $\mathbb{R}^{d}$ case. We consider the
Sobolev space $H^{n}(\mathbb{R}^{d}).$ This space consists of those functions
in $L^{2}(\mathbb{R}^{d})$ all of whose partial derivatives of order at most
$n$ (computed in the distributional sense) are in $L^{2}(\mathbb{R}^{d}).$ Let
$c_{n}$ be a positive constant, whose value will be fixed later. We consider
the inner product on $H^{n}(\mathbb{R}^{d})$ given by%
\[
\left\langle f_{1},f_{2}\right\rangle _{n}=c_{n}\left\langle f_{1}%
,f_{2}\right\rangle _{L^{2}}+(-1)^{n}\left\langle f_{1},\Delta^{n}%
f_{2}\right\rangle _{L^{2}}.
\]
(Technically, this formula assumes that $f_{2}$ is in $H^{2n}(\mathbb{R}%
^{d}),$ but a simple integration-by-parts argument shows that the inner
product extends continuously to $f_{2}$ in $H^{n}(\mathbb{R}^{d}).$) For any
choice of $c_{n}>0,$ $H^{n}(\mathbb{R}^{d})$ is a Hilbert space with respect
to this inner product.

Suppose now that $f$ is \textquotedblleft nice,\textquotedblright\ say, having
Fourier transform that is a $C^{\infty}$ function of compact support. In that
case, all of the integrations by parts that follow are easily justified. Using
the isometry result (Theorem \ref{isometry.thm1}) and the fact that
$\Delta^{n}$ commutes with the heat operator, we obtain%
\[
\left\langle f,\Delta^{n}f\right\rangle _{L^{2}}=\int_{\mathbb{C}^{d}%
}\overline{F(z)}\left(  \Delta_{\mathbb{C}}^{n}F(z)\right)  \frac
{e^{-\left\vert y\right\vert ^{2}/t}}{(\pi t)^{d/2}}~dy~dx,
\]
where%
\[
\Delta_{\mathbb{C}}=\sum_{k=1}^{d}\frac{\partial^{2}}{\partial z_{k}^{2}}.
\]
We now integrate by parts. Since $\overline{F(z)}$ is antiholomorphic, the $z
$-derivatives do not \textquotedblleft see\textquotedblright\ this factor and
all of the derivatives go onto the Gaussian factor, giving%
\[
\left\langle f,\Delta^{n}f\right\rangle _{L^{2}}=\int_{\mathbb{C}^{d}%
}\left\vert F(z)\right\vert ^{2}\left(  \Delta_{\mathbb{C}}^{n}\frac
{e^{-\left\vert y\right\vert ^{2}/t}}{(\pi t)^{d/2}}\right)  ~dy~dx.
\]

Since, now, the Gaussian factor is independent of $x,$ applying $\Delta
_{\mathbb{C}}$ to this factor gives the same result as applying $\Sigma
\lbrack(-i/2)\partial/\partial y_{k}]^{2}=-\frac{1}{4}\Delta_{y}.$ Multiplying
and dividing by the Gaussian factor then gives%
\begin{equation}
\left\langle f,f\right\rangle _{n}=\int_{\mathbb{C}^{d}}\left\vert
F(z)\right\vert ^{2}\left[  c_{n}+h_{n,t}(y)\right]  \frac{e^{-\left\vert
y\right\vert ^{2}/t}}{(\pi t)^{d/2}}~dy~dx,\label{sob.norm}%
\end{equation}
where%
\[
h_{n,t}(y)=\frac{1}{4^{n}}e^{\left\vert y\right\vert ^{2}/t}\Delta
^{n}(e^{-\left\vert y\right\vert ^{2}/t}).
\]
The function $h_{n,t}$ is a polynomial of degree $2n$ whose leading term is a
positive multiple of $\left\vert y\right\vert ^{2n}.$ (Essentially, $h_{n,t}$
is the Wick ordering of the function $\left\vert y\right\vert ^{2n}. $) If we
choose the constant $c_{n}$ large enough, then $c_{n}+h_{n,t}(y)$ will be
strictly positive for all $y.$

Up to now we have been assuming that $f_{1}$ and $f_{2}$ are \textquotedblleft
nice.\textquotedblright\ However, once (\ref{sob.norm}) is established for
nice functions, it is not hard to show that it holds for all functions in
$H^{n}(\mathbb{R}^{d}).$ Furthermore, it can be shown that the map sending $f$
to the analytic continuation of $F:=e^{t\Delta/2}f$ is an isometry of
$H^{n}(\mathbb{R}^{d})$ \textit{onto} the Hilbert space of holomorphic
functions for which (\ref{sob.norm}) is finite. (Compare the reasoning in
Section 10 of \cite{H1}.) Since $h_{n,t}$ behaves like $\left\vert
y\right\vert ^{2n}$ at infinity, the functions satisfying (\ref{sob.norm}) are
precisely the functions satisfying (\ref{sobolev.int1}). This establishes the
first assertion in Theorem \ref{sobolev.thm1}.

To establish the second assertion in Theorem \ref{sobolev.thm1}, we need to
estimate the reproducing kernel for the space of holomorphic functions
satisfying (\ref{sobolev.int1}). This can be done by the argument in Section 2
of \cite{HL}. The details are omitted.
\end{proof}

\begin{proof}[Proof of Theorem \ref{pointwise.thm}] If we write out the formula
(\ref{convolution2}) for $F$ we obtain%
\begin{equation}
F(x+iy)=(2\pi t)^{-d/2}e^{\left\vert y\right\vert ^{2}/2t}e^{-iy\cdot x/t}%
\int_{\mathbb{R}^{d}}e^{iy\cdot x^{\prime}/t}e^{-|x-x^{\prime}|^{2}%
/2t}f(x^{\prime})~dx^{\prime}.\label{f.int}%
\end{equation}
Since $f$ is in $H^{n}(\mathbb{R}^{d})$, it is easily seen that the function
$g_{x}(x^{\prime}):=e^{-|x-x^{\prime}|^{2}/2t}f(x^{\prime})$ is in
$H^{n}(\mathbb{R}^{d})$ for each $x.$ Now, the integral on the right-hand side
of (\ref{f.int}) is (up to a constant) the Fourier transform of $g_{x},$
evaluated at the point $\xi=-y/t.$ Since the Fourier transform $\hat{g}$ of a
function $g$ in $H^{n}(\mathbb{R}^{d})$ satisfies $\int_{\mathbb{R}^{d}%
}\left\vert \hat{g}(\xi)(1+\left\vert \xi\right\vert ^{n})\right\vert
^{2}~d\xi<\infty,$ we conclude that%
\[
\int_{\mathbb{R}^{d}}\left\vert F(x+iy)e^{-\left\vert y\right\vert ^{2}%
/2t}(1+\left\vert y\right\vert ^{n})\right\vert ^{2}~dy<\infty
\]
for all $x.$ The Schwarz inequality then tells us that
\begin{align}
& \int_{\mathbb{R}^{d}}\left\vert F(x+iy)\right\vert e^{-\left\vert
y\right\vert ^{2}/2t}~dy\nonumber\\
& \leq\left[  \int_{\mathbb{R}^{d}}\left(  \frac{1}{1+\left\vert y\right\vert
^{n}}\right)  ^{2}dy\int_{\mathbb{R}^{d}}\left\vert F(x+iy)e^{-\left\vert
y\right\vert ^{2}/2t}(1+\left\vert y\right\vert ^{n})\right\vert
^{2}~dy\right]  ^{1/2},\label{f.ineq}%
\end{align}
where the right-hand side of (\ref{f.ineq}) is finite provided that $2n>d.$
This means that we may take $R=\infty$ in (\ref{inv.int1}) with absolute
convergence of the resulting integral.
\end{proof}

We conclude this section by recording two results of Bargmann, about the
Schwarz space and the space of tempered distributions. Note that tempered
distributions are \textit{less} smooth than $L^{2}$ (they have, roughly, a
negative degree of differentiability) so we get polynomially \textit{worse}
behavior. These results are obtained by adapting Theorem 1.7 of \cite{Ba2} to
our normalization conventions.

\begin{theorem}
[Bargmann]\label{schwarz.thm}A function $F$ on $\mathbb{R}^{d}$ is of the form
$F=e^{t\Delta/2}f,$ with $f$ in the Schwarz space, if and only if for every
$n$ there exists a constant $c_{n}$ with%
\[
\left\vert F(x+iy)\right\vert \leq c_{n}\frac{e^{\left\vert y\right\vert
^{2}/2t}}{(1+\left\vert x\right\vert ^{2}+\left\vert y\right\vert ^{2})^{n}}.
\]
A function $F$ on $\mathbb{R}^{d}$ is of the form $F=e^{t\Delta/2}f,$ with $f$
in the space of tempered distributions, if and only if for some $n$ and some
constant $c$ we have%
\[
\left\vert F(x+iy)\right\vert \leq ce^{\left\vert y\right\vert ^{2}%
/2t}(1+\left\vert x\right\vert ^{2}+\left\vert y\right\vert ^{2})^{n}.
\]

\end{theorem}

\subsection{Multiplication properties}

Our metatheorem says that the functions in the range of the time-$t$ heat
operator are those having an analytic continuation to $\mathbb{C}^{d}$ with
growth in the imaginary directions like $e^{\left\vert y\right\vert ^{2}/2t}.
$ Thus if $F_{1}$ is in the range of the time-$t$ heat operator and $F_{2}$ is
in the range of the time-$s$ heat operator, we expect that the product
function $F_{1}F_{2}$ will have growth in the imaginary directions like
$e^{\left\vert y\right\vert ^{2}/2r},$ where $r$ satisfies
\[
\frac{1}{r}=\frac{1}{s}+\frac{1}{t}.
\]
Thus we expect $F_{1}F_{2}$ to be in the range of the time-$r$ heat operator.
Indeed, a straightforward application of the preceding results yields a
precise result of that sort.

\begin{theorem}
\label{mult.thm1}Suppose that $F_{1}$ is of the form $F_{1}=e^{t\Delta/2}%
f_{1}$ with $f_{1}\in L^{2}(\mathbb{R}^{d},dx)$ and that $F_{2}$ is of the
form $e^{s\Delta/2}f_{2}$ with $f_{2}\in L^{2}(\mathbb{R}^{d},dx).$ Let $r$ be
such that $1/r=1/s+1/t.$ Then there exists a unique $f$ in $L^{2}%
(\mathbb{R}^{d},dx)$ such that%
\[
F_{1}F_{2}=e^{r\Delta/2}f.
\]
Furthermore, if either $f_{1}$ or $f_{2}$ is in the Schwarz space then $f$ is
also in the Schwarz space.
\end{theorem}

\begin{proof}
In light of Theorem \ref{isometry.thm1} (with $t$ replaced by $r$), we need to
verify that the holomorphic function $F_{1}F_{2}$ is square-integrable against
the measure $e^{-\left\vert y\right\vert ^{2}/r}~dy~dx.$ To do this, we use
the pointwise bound (\ref{bound1}) on the function $F_{1}$ together with the
square-integrability property that $F_{2}$ has from Theorem
\ref{isometry.thm1}. This establishes the first assertion in the theorem. The
assertion about Schwarz spaces follows from Theorem \ref{schwarz.thm} and the
pointwise bound (\ref{bound1}).
\end{proof}

\section{The sphere case\label{sphere.sec}}

The main results of this section are special cases of results of Stenzel
\cite{St}, who considered arbitrary compact symmetric spaces. Nevertheless, I
think it is instructive to consider this case separately, since in this case
one does not need the machinery of Lie group theory to state the results. See
\cite{bull} for a survey of related results.

Our general principle is that we expect the functions in the range of the
time-$t$ heat operator to have analytic continuations with the same rate of
growth in the imaginary directions as the time-$t$ heat kernel. Meanwhile, the
Minakshisundaram--Pleijel asymptotic series for the heat kernel on a manifold
suggests that the heat kernel should behave roughly like a Gaussian times
$j^{-1/2},$ where $j$ is the Jacobian of the exponential mapping.
Specifically, the zeroth term in the series is precisely of this form. If we
are on a nice manifold such as a sphere, we may hope that this leading-order
approximation will be accurate even after analytic continuation. Thus, in the
sphere case, the growth rate in our metatheorem will come directly from the
analytic continuation of a Gaussian times $j^{-1/2}.$ The Jacobian factor
$j^{-1/2}$ contributes a vitally important exponentially decaying factor to
our metatheorem in the sphere case, a factor that we do not have in the
Euclidean case.

Before stating our metatheorem precisely, we need a bit of set-up. First, we
need the appropriate manifold to which functions in the range of the heat
operator for $S^{d}$ should be analytically continued. This will be the
\textquotedblleft complexified sphere\textquotedblright\ $S_{\mathbb{C}}^{d}.
$ If the ordinary sphere $S^{d}$ is
\[
S^{d}=\left\{  \left.  x\in\mathbb{R}^{d+1}\right\vert x_{1}^{2}%
+\cdots+x_{d+1}^{2}=1\right\}  ,
\]
then we define the complexified sphere $S_{\mathbb{C}}^{d}$ to be%
\[
S_{\mathbb{C}}^{d}=\left\{  \left.  z\in\mathbb{C}^{d+1}\right\vert z_{1}%
^{2}+\cdots+z_{d+1}^{2}=1\right\}  .
\]
Note that there are no complex-conjugates in the definition of $S_{\mathbb{C}%
}^{d}$; it is \textit{not} the unit sphere in $\mathbb{C}^{d+1}\cong
\mathbb{R}^{2(d+1)}.$ Rather, $S_{\mathbb{C}}^{d}$ is a noncompact complex
manifold with complex dimension $d,$ defined by a single holomorphic condition
in $\mathbb{C}^{d+1}.$ The real sphere $S^{d}$ sits inside $S_{\mathbb{C}}%
^{d}$ as a totally real submanifold of maximal dimension (i.e., it looks
locally like $\mathbb{R}^{d}$ inside $\mathbb{C}^{d}$). Thus it makes sense to
speak of analytically continuing (sufficiently regular) functions from $S^{d}$
to $S_{\mathbb{C}}^{d}.$

Second, we need the right way to measure the rate of growth of functions on
$S_{\mathbb{C}}^{d}.$ This comes by identifying $S_{\mathbb{C}}^{d}$ with the
tangent bundle $T(S^{d}).$ We think of the tangent bundle to $S^{d}$ as the
set%
\[
T(S^{d})=\left\{  \left.  (x,Y)\in\mathbb{R}^{d+1}\times\mathbb{R}%
^{d+1}\right\vert ~\left\vert x\right\vert =1,~Y\cdot x=0\right\}  .
\]
Then we have the \textquotedblleft exponential map\textquotedblright\ for
$S^{d},$ in the differential geometric sense. Given a point $x$ in $S^{d}$ and
a tangent vector $Y$ at $x,$ let $\gamma:\mathbb{R}\rightarrow S^{d}$ be the
unique geodesic with $\gamma(0)=x$ and $\left.  d\gamma(\tau)/d\tau\right\vert
_{\tau=0}=Y.$ Then $\exp_{x}(Y)$ is defined to be $\gamma(1).$ In the case of
the sphere, the geodesics are simply great circles and the exponential map may
be computed explicitly as%
\[
\exp_{x}(Y)=\cos\left\vert Y\right\vert ~x+\frac{\sin\left\vert Y\right\vert
}{\left\vert Y\right\vert }Y.
\]

For each $x\in S^{d},$ $\exp_{x}$ is a smooth map of $T_{x}(S^{d})$ into
$S^{d}.$ For each fixed $x,$ we may analytically continue $\exp_{x}$ into a
holomorphic map of the complexified tangent space $T_{x}(S^{d})_{\mathbb{C}}$
into the complexified sphere $S_{\mathbb{C}}^{d},$ where $T_{x}(S^{d}%
)_{\mathbb{C}}=\left\{  \left.  Y\in\mathbb{C}^{d+1}\right\vert Y\cdot
x=0\right\}  .$ To see this, note that both $\cos\theta$ and $\sin
\theta/\theta$ are even functions of $\theta.$ Thus, $\cos\left\vert
Y\right\vert $ and $\sin\left\vert Y\right\vert /\left\vert Y\right\vert $ can
be expressed as power series (with infinite radius of convergence) in powers
of $\left\vert Y\right\vert ^{2}=Y_{1}^{2}+\cdots+Y_{d+1}^{2}.$ For $Y$ in the
complexified tangent space at $x,$ we compute the complex-valued quantity
$Y_{1}^{2}+\cdots+Y_{d+1}^{2}$ and plug this into the same power series, in
order to compute $\exp_{x}(Y)$ as an element of $\mathbb{C}^{d+1}.$ Since
$\exp_{x}$ maps the real tangent space into $S^{d}\subset\mathbb{R}^{d+1},$
its analytic continuation will map the complexified tangent space into
$S_{\mathbb{C}}^{d}\subset\mathbb{C}^{d+1}.$

We then define a map $\Phi:T(S^{d})\rightarrow S_{\mathbb{C}}^{d}$ by the
formula%
\[
\Phi(x,Y)=\exp_{x}(iY).
\]
This may be computed explicitly as%
\[
\Phi(x,Y)=\cosh\left\vert Y\right\vert ~x+i\frac{\sinh\left\vert Y\right\vert
}{\left\vert Y\right\vert }Y
\]
It is not hard to show that $\Phi$ is a diffeomorphism of $T(S^{d})$ with
$S_{\mathbb{C}}^{d}.$ (Compare Section 2 of \cite{St}.) This identification of
$T(S^{d})$ with $S_{\mathbb{C}}^{d}$ can be obtained naturally from the theory
of \textquotedblleft adapted complex structures\textquotedblright\ as
developed in \cite{GS1,GS2,LS,Sz}. (See also Section 3 of \cite{H3}.)

We are now ready to state our metatheorem for the $S^{d}$ case.

\begin{metatheorem}
\label{sphere.meta}The functions in the range of the heat operator
$e^{t\Delta/2}$ for $S^{d}$ are those functions $F$ having an analytic
continuation to $S_{\mathbb{C}}^{d}$ such that $F(\exp_{x}(iY))$ grows at most
like%
\begin{equation}
e^{\left\vert Y\right\vert ^{2}/2t}\left(  \frac{\left\vert Y\right\vert
}{\sinh\left\vert Y\right\vert }\right)  ^{(d-1)/2}.\label{sd.metabound}%
\end{equation}

\end{metatheorem}

As in the $\mathbb{R}^{d}$ case, we will allow ourselves some flexibility in
how we interpret the metatheorem, sometimes in an $L^{2}$ sense and sometimes
in a pointwise sense, and always allowing polynomial variations depending on
the \textit{domain} of the heat operator ($L^{2},$ $C^{\infty}, $ etc.).

Note that (except when $d=1$) the function $(\left\vert Y\right\vert
/\sinh\left\vert Y\right\vert )^{(d-1)/2}$ has exponential decay as
$\left\vert Y\right\vert $ goes to infinity. Although the exponential decay of
this factor in (\ref{sd.metabound}) might seem insignificant compared to the
Gaussian growth of the factor $e^{\left\vert Y\right\vert ^{2}/2t},$ the
exponential factor is in fact very important. After all, we expect, based on
our experience in the $\mathbb{R}^{d}$ case, that smoothness properties of the
initial data $f$ will reflect themselves in polynomial fluctuations around the
\textquotedblleft basic\textquotedblright\ rate of growth for $F.$ If we miss
an exponential factor in the basic growth rate, we are not going to be able to
see the polynomial fluctuations.

The factor $(\left\vert Y\right\vert /\sinh\left\vert Y\right\vert
)^{(d-1)/2}$ is nothing but the analytic continuation of the Jacobian factor
$j^{-1/2}$ referred to at the beginning of this section. Specifically, the
Jacobian of the exponential mapping at a point $x$ in $S^{d}$ is given by
\[
j(\exp_{x}Y)=\left(  \frac{\sin\left\vert Y\right\vert }{\left\vert
Y\right\vert }\right)  ^{d-1}.
\]
Thus the extra factor in our metatheorem is nothing but $j^{-1/2}(\exp
_{x}(iY)).$ We expect (based on the Minakshisundaram--Pleijel expansion) the
heat kernel $\rho_{t,x}$ at a point $x$ in $S^{d}$ to satisfy
\[
\rho_{t,x}(\exp_{x}Y)\approx(2\pi t)^{-d/2}e^{-\left\vert Y\right\vert
^{2}/2t}j^{-1/2}(\exp_{x}Y).
\]
Thus we expect the analytic continuation of $\rho_{t}$ to grow in the
imaginary directions like (a constant times) $e^{\left\vert Y\right\vert
^{2}/2t}j^{-1/2}(\exp_{x}(iY))$. This line of reasoning accounts for the form
of our metatheorem.

\subsection{Isometry and inversion formulas}

Consider the spherical Laplacian $\Delta$ on $S^{d},$ as described, for
example, in \cite{Ta}. This is just the Laplace--Beltrami operator for $S^{d}
$ as a Riemannian manifold, except that we choose the sign so that $\Delta$ is
a negative operator. This operator is essentially self-adjoint as an unbounded
operator on $L^{2}(S^{d})$ with domain $C^{\infty}(S^{d}).$ We let $\Delta$
also denote the unique self-adjoint extension of this operator. The spectral
theorem allows us define $e^{t\Delta/2}$ as a contraction operator on
$L^{2}(S^{d}).$ For each point $x\in S^{d},$ there is the \textquotedblleft
heat kernel based at $x$\textquotedblright, denoted $\rho_{t,x}$ with the
property that%
\[
(e^{t\Delta/2}f)(x)=\int_{S^{d}}\rho_{t,x}(y)f(y)~dy.
\]

Note that under the diffeomorphism $\Phi$ of $T(S^{d})$ with $S_{\mathbb{C}%
}^{d},$ the fiber directions in $T(S^{d})$ correspond to what we may call the
\textquotedblleft imaginary directions\textquotedblright\ in $S_{\mathbb{C}%
}^{d},$ namely, the directions corresponding to the exponential map with a
pure imaginary argument. The isometry and inversion formula involve an
identification of fiber directions in $T(S^{d}),$ and thus of the imaginary
directions in $S_{\mathbb{C}}^{d},$ with $d$-dimensional hyperbolic space. The
\textit{heat kernel measure} for hyperbolic space will then play the role of
the Gaussian measure in the imaginary direction that we have in the
$\mathbb{R}^{d}$ case.

Concretely, this means that we consider the function $\nu_{t}(R)$ satisfying
the differential equation%
\begin{equation}
\frac{\partial}{\partial t}\nu_{t}(R)=\frac{1}{2}\left[  \frac{\partial^{2}%
}{\partial R^{2}}+(d-1)\frac{\cosh R}{\sinh R}\frac{\partial}{\partial
R}\right]  \nu_{t}(R)\label{nut.def}%
\end{equation}
subject to the initial condition%
\begin{equation}
\lim_{t\searrow0}c_{d}\int_{0}^{\infty}f(R)\nu_{t}(R)(\sinh R)^{d-1}%
dR=f(0)\label{nut.init}%
\end{equation}
for all continuous functions $f$ on $[0,\infty)$ with at most exponential
growth. Here $c_{d}$ is the volume of the unit sphere in $\mathbb{R}^{d}.$ The
function $\nu_{t}(R)$ is just the heat kernel function on $d$-dimensional
hyperbolic space, as a function of the distance $R$ from the basepoint. The
operator in square brackets in (\ref{nut.def}) is just the hyperbolic
Laplacian acting on radial functions. The integral of $f(R)$ against
$c_{d}(\sinh R)^{d-1}dR$ in (\ref{nut.init}) is the formula for the integral
of a radial function on hyperbolic space. See \cite[Sect. 5.7]{Da} or
\cite[Sect. VI]{HM1} for more information and for formulas for $\nu_{t}$ in
various dimensions.

We are now ready to state the isometry and inversion formulas for the sphere
case. These results are special cases of results of Stenzel \cite{St}, who
considers general compact symmetric spaces. The case of $S^{3}$ (which is
isometric to the group $SU(2)$ with a bi-invariant Riemannian metric) falls
under the earlier work of \cite{H1,H2} on the compact group case. A
self-contained treatment of the isometry formula in the sphere case was given
(subsequently to \cite{St}) in \cite{HM1}. See also \cite{KR2,HM2} for more on
the sphere case and \cite[Sect. IV.D]{Th} for the special case of $S^{3}\cong
SU(2).$

\begin{theorem}
[Isometry]A function $F$ on $S^{d}$ is of the form $F=e^{t\Delta/2}f,$ with
$f\in L^{2}(S^{d}),$ if and only if $F$ has an analytic continuation to
$S_{\mathbb{C}}^{d}$ satisfying%
\begin{equation}
\int_{x\in S^{d}}\int_{Y\in T_{x}(S^{d})}\left\vert F\left(  \exp_{x}\left(
iY\right)  \right)  \right\vert ^{2}\nu_{2t}(\left\vert 2Y\right\vert
)~\left(  \frac{\sinh\left\vert 2Y\right\vert }{\left\vert 2Y\right\vert
}\right)  ^{d-1}2^{d}dY~dx<\infty.\label{isometry.intsd}%
\end{equation}
Here $dx$ denotes the natural volume measure on $S^{d}$ and $dY$ denotes
Lebesgue measure on $T_{x}(S^{d}).$ If $F$ is such a function then the
integral in (\ref{isometry.intsd}) is equal to $\int_{S^{d}}|f(x)|^{2}~dx.$
\end{theorem}

\begin{theorem}
[Inversion]For any $f$ in $L^{2}(S^{d}),$ let $F=e^{t\Delta/2}f.$ Then%
\[
f(x)=\lim_{R\rightarrow\infty}\int_{\substack{Y\in T_{x}(S^{d}) \\\left\vert
Y\right\vert \leq R}}F\left(  \exp_{x}(iY)\right)  \nu_{t}(\left\vert
Y\right\vert )\left(  \frac{\sinh\left\vert Y\right\vert }{\left\vert
Y\right\vert }\right)  ^{d-1}dY,
\]
where the limit is in the norm topology of $L^{2}(S^{d})$ and where $dY$
denotes Lebesgue measure on $T_{x}(S^{d}).$ Furthermore, if $f$ is
sufficiently smooth then%
\[
f(x)=\int_{Y\in T_{x}(S^{d})}F\left(  \exp_{x}(iY)\right)  \nu_{t}(\left\vert
Y\right\vert )\left(  \frac{\sinh\left\vert Y\right\vert }{\left\vert
Y\right\vert }\right)  ^{d-1}dY
\]
with the integral on the right-hand side being absolutely convergent for all
$x\in S^{d}.$
\end{theorem}

These theorems take on an especially explicit form in the case $d=3$; see
Section \ref{s3.sect}. The expression $(\sinh\left\vert Y\right\vert
/\left\vert Y\right\vert )^{d-1}~dY$ is the formula for the Riemannian volume
measure on hyperbolic space, written in global exponential coordinates. Thus
the quantity $\nu_{t}(\left\vert Y\right\vert )(\sinh\left\vert Y\right\vert
/\left\vert Y\right\vert )^{d-1}~dY$ is the \textit{heat kernel measure} for
hyperbolic space (i.e., the heat kernel function times the Riemannian volume measure).

The idea behind the inversion formula is roughly that the integration against
the heat kernel is computing the forward hyperbolic heat operator in the
imaginary variables. A fairly simple analytic continuation argument shows that
on holomorphic functions, the \textit{forward} \textit{hyperbolic} heat
operator in the imaginary directions is the same as the \textit{backward}
\textit{spherical} heat operator in the real directions. Thus the integration
in the inversion formula undoes the heat equation that takes us from $f$ to
$F.$ This line of reasoning should be compared to the heuristic argument for
the inversion formula (Theorem \ref{inv.thm1}) in the $\mathbb{R}^{d}$ case.
Stenzel makes this argument rigorous and then uses the inversion formula to
obtain the isometry formula (as in \cite{H2} in the compact group case).

Observe that in the inversion formula we have $\nu_{t}(\left\vert Y\right\vert
),$ whereas in the isometry formula we have $\nu_{2t}(\left\vert 2Y\right\vert
).$ In \cite{St}, the inversion formula is obtained first and the isometry
formula is reduced to the inversion formula. In the process of this reduction,
the scalings $Y\rightarrow2Y$ and $t\rightarrow2t$ arise naturally. It is hard
to predict the correct way to scale the variables by looking at the Euclidean
case, since in that case a scaling of the space variables can be absorbed into
a scaling of the time variables. That is, in the $\mathbb{R}^{d}$ case we have
$2^{d}\nu_{2t}(2y)=\nu_{t/2}(y),$ which definitely does not hold in the
hyperbolic case.

Let us think about how the isometry formula fits with our metatheorem in the
sphere case. It is known \cite{DM, HSt} that the heat kernel $\nu_{t}$
satisfies, for each $t,$
\begin{equation}
\nu_{t}(R)\approx e^{-R^{2}/2t}\left(  \frac{R}{\sinh R}\right)
^{(d-1)/2}.\label{nut.approx1}%
\end{equation}
Thus
\begin{equation}
\nu_{2t}(\left\vert 2Y\right\vert )\left(  \frac{\sinh\left\vert 2Y\right\vert
}{\left\vert 2Y\right\vert }\right)  ^{d-1}\approx e^{-\left\vert Y\right\vert
^{2}/t}\left(  \frac{\sinh\left\vert 2Y\right\vert }{\left\vert 2Y\right\vert
}\right)  ^{(d-1)/2}.\label{nut.approx2}%
\end{equation}
Here $\approx$ means that the ratio of the two sides is bounded and bounded
away from zero as a function of $R$ for each fixed $t.$ Thus if $\left\vert
F\right\vert ^{2}$ is to be \textit{square}-integrable against the measure in
(\ref{isometry.intsd}) we must have, roughly,
\begin{equation}
\left\vert F(\exp_{x}(iY))\right\vert \lesssim e^{\left\vert Y\right\vert
^{2}/2t}\left(  \frac{\left\vert 2Y\right\vert }{\sinh\left\vert 2Y\right\vert
}\right)  ^{(d-1)/4}.\label{approx.bnd}%
\end{equation}
Since $(\sinh\left\vert 2Y\right\vert )^{(d-1)/4}$ has the same behavior at
infinity as $(\sinh\left\vert Y\right\vert )^{(d-1)/2},$ the right-hand side
of (\ref{approx.bnd}) has the same behavior as the quantity in our
metatheorem, to within a polynomial factor in $\left\vert Y\right\vert .$

Let us think, in turn, about how the inversion formula fits with our
metatheorem. To get convergence we need something like
\[
\left\vert F(\exp_{x}(iY))\right\vert \lesssim\left[  \nu_{t}(\left\vert
Y\right\vert )\left(  \frac{\sinh\left\vert Y\right\vert }{\left\vert
Y\right\vert }\right)  ^{d-1}\right]  ^{-1}.
\]
Putting in the estimate (\ref{nut.approx1}), this comes to%
\[
\left\vert F(\exp_{x}iY)\right\vert \lesssim e^{\left\vert Y\right\vert
^{2}/2t}\left(  \frac{\left\vert Y\right\vert }{\sinh\left\vert Y\right\vert
}\right)  ^{(d-1)/2},
\]
which is precisely what we have in the metatheorem.

\subsection{Pointwise bounds}

The optimal pointwise bounds for holomorphic functions $F$ on $S_{\mathbb{C}%
}^{d}$ satisfying (\ref{isometry.intsd}) are obtained from the reproducing
kernel for the corresponding $L^{2}$ space of holomorphic functions. It
follows easily from the results of Stenzel (using the method in Section 9 of
\cite{H1}) that the reproducing kernel for the space of holomorphic functions
satisfying (\ref{isometry.intsd}) is given in terms of the analytic
continuation of the heat kernel $\rho_{t}$ on $S^{d}.$ This leads to the following.

\begin{theorem}
\label{sdbnd.thm}For $f\in L^{2}(S^{d}),$ let $F=e^{t\Delta/2}f.$ Then the
analytic continuation of $F$ satisfies%
\[
\left\vert F(\exp_{x}(iY)\right\vert \leq\left\Vert f\right\Vert _{L^{2}%
(S^{d})}\sqrt{\rho_{2t,x}(\exp_{x}(2iY))}.
\]

\end{theorem}

I conjecture that for each $d$ and $t,$ there exists a constant $a_{t,d}$ such
that%
\begin{equation}
\rho_{t,x}(\exp_{x}(iY))\leq a_{t,d}e^{\left\vert Y\right\vert ^{2}/2t}\left(
\frac{\left\vert Y\right\vert }{\sinh\left\vert Y\right\vert }\right)
^{(d-1)/2}.\label{rhot.bnd}%
\end{equation}
Such bounds are known for the cases $d=1$ and $d=3.$ (Apply Theorem 2 of
\cite{H3} in the case $K=S^{1}$ or $K=SU(2)\cong S^{3}$.) If this conjecture
were true, we would obtain the following precise pointwise version of the
estimate (\ref{approx.bnd}):%
\begin{equation}
\left\vert F(\exp_{x}(iY))\right\vert \leq\left\Vert f\right\Vert
_{L^{2}(S^{d})}\sqrt{a_{2t,d}}~e^{\left\vert Y\right\vert ^{2}/2t}\left(
\frac{\left\vert 2Y\right\vert }{\sinh\left\vert 2Y\right\vert }\right)
^{(d-1)/4}.\label{precise.bnd}%
\end{equation}
Note that this bound is better by a polynomial factor than the one in our
metatheorem, where the bound in our metatheorem is nothing but what one would
expect from the behavior of $\rho_{t}$ in (\ref{rhot.bnd}). This improvement
presumably reflects a slight smoothing that takes place when integrating a
function $f$ in $L^{2}(S^{d})$ against the heat kernel: The behavior of
$e^{t\Delta/2}f$, for $f\in L^{2}(S^{d})$, is slightly better than the
behavior of the heat kernel $\rho_{t}.$

We also expect, based on our experience in the $\mathbb{R}^{d}$ case, that
assuming smoothness for $f$ gives polynomial improvement in the pointwise
behavior of $F$ (as compared to (\ref{precise.bnd})). This, again, is known in
the cases $d=1$ and $d=3,$ by applying the results of \cite{HL} to the case
$K=S^{1}$ or $K=SU(2)\cong S^{3}$. See Section \ref{s3.sect} for more on the
$S^{3}$ case.

\subsection{Multiplication properties}

The metatheorem in the sphere case suggests a multiplication result similar to
what we have in the Euclidean case (Theorem \ref{mult.thm1}). However, the
\textquotedblleft non-Euclidean\textquotedblright\ factor involving
$\sinh\left\vert Y\right\vert $ works to our advantage in this case. If we
multiply to functions in the range of the heat operator, the Gaussian factors
combine in the same way as in the $\mathbb{R}^{d}$ case, but we now get the
exponentially decaying sinh factors twice. This suggests the following result,
which can be proved by lifting to the group $SO(d+1)$ and using results of
\cite{HL}. (See the proof of Theorem \ref{mult.thm3} in Section
\ref{compact.sec}.)

\begin{theorem}
\label{mult.thm2}Assume $d\geq2.$ Given $f_{1},f_{2}$ in $L^{2}(S^{d}),$ let
$F_{1}=e^{t\Delta/2}f_{1}$ and $F_{2}=e^{s\Delta/2}f_{2},$ for two positive
numbers $s$ and $t.$ Let $r$ be such that $1/r=1/s+1/t.$ Then there exists a
unique $f$ in $C^{\infty}(S^{d})$ such that%
\[
e^{r\Delta/2}f=F_{1}F_{2}.
\]

\end{theorem}

Note that although $f_{1}$ and $f_{2}$ are assumed only to be in $L^{2},$ we
obtain that $f$ is in $C^{\infty}.$

Although this theorem has the look of the sort of soft result that might hold
quite generally, it is actually very special. For example, nothing remotely
like this holds in the hyperbolic case. Indeed, the class of manifolds for
which a result of this type holds is probably very small.

\subsection{The $S^{3}$ case\label{s3.sect}}

In the case of the 3-sphere, the isometry and inversion formulas take on a
particularly explicit form. Furthermore, in this case, we can obtain results
involving smoothness of the initial data similar to what we have in the
$\mathbb{R}^{d}$ case. Such results are expected to hold for other spheres,
but the proofs are more difficult when $d\neq3,$ due to the less explicit
formulas for the relevant heat kernels. What is special about the $S^{3}$ case
is that $S^{3}$ is a group. That is, $S^{3}$ with the standard metric is
isometric to the group $SU(2)$ with a bi-invariant metric. The results of this
subsection come from \cite{H1,H2,H3, HL}.

\begin{theorem}
\label{s3.thm}In the $S^{3}$ case, the isometry and inversion formulas take
the form%
\begin{equation}
\int_{S^{d}}\left\vert f(x)\right\vert ^{2}dx=\int_{x\in S^{3}}\int_{Y\in
T_{x}(S^{3})}\left\vert F\left(  \exp_{x}\left(  iY\right)  \right)
\right\vert ^{2}e^{-t}\frac{\sinh\left\vert 2Y\right\vert }{\left\vert
2Y\right\vert }\frac{e^{-\left\vert Y\right\vert ^{2}/t}}{(\pi t)^{3/2}%
}~dY~dx\label{isometry.s3}%
\end{equation}
and%
\begin{equation}
f(x)=\int_{Y\in T_{x}(S^{3})}F(\exp_{x}(iY))e^{-t/2}\frac{\sinh\left\vert
Y\right\vert }{\left\vert Y\right\vert }\frac{e^{-\left\vert Y\right\vert
^{2}/2t}}{(2\pi t)^{3/2}}~dY.\label{inversion.s3}%
\end{equation}
Here, as usual, $F$ is the analytic continuation of $e^{t\Delta/2}f.$ The
isometry formula (\ref{isometry.s3}) holds for all $f$ in $L^{2}(S^{d}),$ and
$F$ is in the image of $L^{2}(S^{d})$ under $e^{t\Delta/2}$ if and only if the
integral on the right-hand side of (\ref{isometry.s3}) is finite. The
inversion formula (\ref{inversion.s3}) holds for sufficiently smooth $f$, with
absolute convergence of the integral for each $x.$
\end{theorem}

In addition, the expected results (\ref{rhot.bnd}) and (\ref{precise.bnd}) are
known to hold in the $d=3$ case.

By specializing the results of \cite{HL} to the case $K=SU(2),$ we obtain the
following counterpart to Theorem \ref{sobolev.thm1}.

\begin{theorem}
\label{sobolev.thm2}The function $f$ is in the Sobolev space $H^{2n}(S^{3})$
if and only if $F$ satisfies%
\begin{equation}
\int_{x\in S^{3}}\int_{Y\in T_{x}(S^{3})}\left\vert F\left(  \exp_{x}\left(
iY\right)  \right)  \right\vert ^{2}(1+\left\vert Y\right\vert ^{4n}%
)e^{-t}\frac{\sinh\left\vert 2Y\right\vert }{\left\vert 2Y\right\vert }%
\frac{e^{-\left\vert Y\right\vert ^{2}/t}}{(\pi t)^{3/2}}~dY~dx<\infty
.\label{s3.hn}%
\end{equation}
For such an $f$ we have the pointwise bounds%
\begin{equation}
\left\vert F(\exp_{x}(iY))\right\vert \leq Ce^{\left\vert Y\right\vert
^{2}/2t}\left(  \frac{\left\vert 2Y\right\vert }{\sinh\left\vert 2Y\right\vert
}\right)  ^{1/2}\frac{1}{1+\left\vert Y\right\vert ^{2n}}.\label{s3.hnbnd}%
\end{equation}

A function $f$ is in $C^{\infty}(S^{3})$ if and only if $F$ satisfies a bound
of the form (\ref{s3.hnbnd}) for each $n=1,2,\ldots$ (with the constant $C$
depending on $n$).
\end{theorem}

\section{Compact symmetric spaces\label{compact.sec}}

The concerning spheres in the previous section have straightforward extensions
to arbitrary compact symmetric spaces. I briefly summarize those extensions
here. For more details, see \cite{St} or Section 5 of \cite{HM3}. I will for
simplicity restrict to the case of simply connected compact symmetric spaces.
These spaces may be represented in the form $X=U/K,$ where $U$ is a simply
connected compact Lie group (automatically semisimple), $K$ is the fixed-point
subgroup of an involution of $U,$ and the metric on $X$ is invariant under the
action of $U.$ Conversely, any Riemannian manifold constructed in this way is
a simply connected compact symmetric space. (See \cite{He} for details.) We
will assume, without loss of generality, that the action of $U$ on $X$ is
locally effective, that is, that the group of elements acting trivially on $X$
is discrete. In that case, $U$ and $K$ are unique up to isomorphism (for a
given $X$) and $U$ is the universal cover of the identity component of the
isometry group of $X.$

Let $\mathfrak{u}$ denote the Lie algebra of $U$ and let $\mathfrak{u}%
_{\mathbb{C}}:=\mathfrak{u}+i\mathfrak{u}$ be the complexification of
$\mathfrak{u}.$ Let $U_{\mathbb{C}}$ be the unique simply connected Lie group
with Lie algebra $\mathfrak{u}_{\mathbb{C}}.$ Let $K_{\mathbb{C}}$ be the
connected Lie subgroup of $U_{\mathbb{C}}$ whose Lie algebra is $\mathfrak{k}%
_{\mathbb{C}}:=\mathfrak{k}+i\mathfrak{k},$ where $\mathfrak{k}$ is the Lie
algebra of $K.$ Then we define the \textquotedblleft
complexification\textquotedblright\ $X_{\mathbb{C}}$ of $X$ to be the complex
manifold $U_{\mathbb{C}}/K_{\mathbb{C}}.$ As in the sphere case, we have a
diffeomorphism $\Phi$ of $T(U/K)$ with $U_{\mathbb{C}}/K_{\mathbb{C}} $ given
by%
\[
\Phi(x,Y)=\exp_{x}(iY),
\]
where $\exp_{x}(iY)$ refers to the analytically continued exponential map for
$U/K.$

The involution of $U$ induces a involution of $\mathfrak{u},$ and
$\mathfrak{u}$ then decomposes as $\mathfrak{u}=\mathfrak{k}+\mathfrak{p},$
where $\mathfrak{p}$ is the $-1$ eigenspace for the involution. If we set
$\mathfrak{g}=\mathfrak{k}+i\mathfrak{p},$ then $\mathfrak{g}$ is a subalgebra
of $\mathfrak{u}.$ Let $G$ be the connected Lie subgroup of $U_{\mathbb{C}}$
whose Lie algebra is $\mathfrak{g}.$ Then $G/K$ (with an appropriately
normalized $G$-invariant Riemannian metric) is a noncompact symmetric space,
called the \textquotedblleft dual\textquotedblright\ of $U/K.$ For example, if
$U/K$ is a $d$-sphere then $G/K$ will be hyperbolic $d$-space.

Let $x_{0}$ denote the identity coset in $U/K\subset U_{\mathbb{C}%
}/K_{\mathbb{C}}$ and let us identify the tangent space $T_{x_{0}}(U/K)$ with
$\mathfrak{p}.$ Since the geometric and group-theoretic exponential mappings
coincide for symmetric spaces, we have%
\[
\exp_{x_{0}}(iY)=e^{iY}\cdot x_{0}%
\]
for all $Y\in\mathfrak{p},$ where $e^{iY}$ is the group-theoretic exponential
of $iY$ inside $U_{\mathbb{C}}.$ Note that $iY\in i\mathfrak{p}\subset
\mathfrak{g}$ and therefore $e^{iY}\in G.$ From this it is not hard to see
that the image of $T_{x_{0}}(U/K)$ under $\Phi$ is precisely the $G$-orbit of
$x_{0}.$ But the stabilizer of $x_{0}$ inside $G$ is just $K,$ and so the
image of $T_{x_{0}}(U/K)$ is naturally identified with the noncompact
symmetric space $G/K.$ This identification is actually just the geometric
exponential mapping for $G/K,$ viewing $\mathfrak{p}\cong i\mathfrak{p}$ as
the tangent space at the identity coset to $G/K.$

Of course, the tangent space at any point $x$ in $U/K$ may be identified
(non\allowbreak uniquely) with $T_{x_{0}}(U/K)$ by means of the action of $U.
$ Thus \textit{every} fiber in the tangent bundle of $U/K$ may be identified
with the noncompact symmetric space $G/K.$ This identification is not unique,
but is unique up to the action of $K$ on $G/K.$ Once each fiber in $T(U/K)$ is
identified with $G/K,$ we may introduce on each fiber the \textit{heat kernel
function }$\nu_{t}$ for $G/K$ and the \textit{Jacobian }$j$\textit{\ of the
exponential mapping} for $G/K.$

We are now ready to state the isometry and inversion formulas for $U/K.$

\begin{theorem}
[Stenzel]1. The \textbf{isometry formula}. A function $F$ on $U/K$ is of the
form $F=e^{t\Delta/2}f$, with $f\in L^{2}(U/K),$ if and only if $F$ has an
analytic continuation to $U_{\mathbb{C}}/K_{\mathbb{C}}$ satisfying%
\[
\int_{x\in U/K}\int_{Y\in T_{x}(U/K)}\left\vert F(\exp_{x}(iY))\right\vert
^{2}\nu_{2t}(2Y)j(2Y)~2^{d}dY~dx.
\]
If $F$ is such a function, then the above integral is equal to $\int
_{U/K}\left\vert f(x)\right\vert ^{2}~dx.$ Here $d=\dim(U/K),$ $dY$ is the
Lebesgue measure on $T_{x}(U/K)$ and $dx$ is the Riemannian volume measure on
$U/K.$

2. The \textbf{inversion formula}. If $f$ is sufficiently smooth and
$F:=e^{t\Delta/2}f,$ then%
\[
f(x)=\int_{T_{x}(U/K)}F(\exp_{x}(iY))\nu_{t}(Y)j(Y)~dY,
\]
with absolute convergence of the integral for all $x$ in $U/K.$
\end{theorem}

Note that in the isometry formula we have $\nu_{2t}(2Y)j(2Y),$ whereas in the
inversion formula we have $\nu_{t}(Y)j(Y).$ Stenzel first proves the inversion
formula and then (as in \cite{H2} in the compact group case) derives the
isometry formula from the inversion formula. In this derivation, the scalings
by a factor of 2 in the space and time variables arise naturally.

Among compact symmetric spaces are compact Lie groups with bi-invariant
metrics. The compact group case, which was considered prior to \cite{St} in
\cite{H1,H2}, is special in various ways. In the first place, there are simple
explicit formulas for $\nu_{t}$ (and $j$) in this case. This allows for
precise estimates of various quantities of interest \cite{H3,HL} that are not
yet known in general. Note that the only spheres that fall into the group case
are $S^{1}$ and $S^{3}.$ This accounts for the particularly explicit nature of
the formulas for the $S^{d}$ case when $d=3.$ In the second place, the group
case is special because in this case, the isometry and inversion formulas can
be obtained by means of geometric quantization \cite{geoquant}. (See also
\cite{FMMN,FMMN2}.) The results of \cite{geoquant} are simply false as soon as
one moves outside of the group case.

Of course, one can always lift analysis on $U/K$ up to the compact group $U.$
This strategy works better for some problems than for others. The
multiplication properties for the range of the heat operator is one problem
where this lifting strategy works well. We obtain, then, the following result.

\begin{theorem}
\label{mult.thm3} Let $U/K$ be a simply connected compact symmetric space.
Suppose that $F_{1}$ is of the form $F_{1}=e^{t\Delta/2}f_{1}$ with $f_{1}\in
L^{2}(U/K)$ and that $F_{2}$ is of the form $e^{s\Delta/2}f_{2}$ with
$f_{2}\in L^{2}(U/K).$ Let $r$ be such that $1/r=1/s+1/t.$ Then there exists a
unique $f$ in $C^{\infty}(U/K)$ such that%
\[
F_{1}F_{2}=e^{r\Delta/2}f.
\]

\end{theorem}

In the nonsimply connected case, the same result holds, \textit{provided} that
one assumes $U/K$ is of the \textquotedblleft compact type\textquotedblright%
\ in the sense of Helgason. This last stipulation is the reason for the
assumption $d\geq2$ in Theorem \ref{mult.thm2}. Without the compact-type
assumption, one gets only the weaker form of the multiplication theorem one
has in the $\mathbb{R}^{d}$ case.

\begin{proof}
We regard $f_{1}$ and $f_{2}$ as right-$K$-invariant functions on $U$ and
apply the pointwise bounds of \cite[Thm. 2]{H3} to $F_{1}$ and $F_{2}.$ Then
the characterization in \cite[Thm. 5]{HL} of the image under $e^{r\Delta/2}$
of $C^{\infty}(U)$ shows that $F_{1}F_{2}=e^{r\Delta/2}f$ for some $f$ in
$C^{\infty}(U).$ Finally, $F_{1}F_{2}$ is right-$K$-invariant and so,
therefore, is $f.$ Thus $f$ may be thought of as a $C^{\infty}$ function on
$U/K.$
\end{proof}

I conclude this section with a conjecture involving the heat kernel on a
compact symmetric space and the heat kernel on the dual noncompact symmetric
space. The conjecture is established in \cite{H3} in the compact group case,
where it is interpreted as a uniform bound on \textquotedblleft phase space
probability densities,\textquotedblright\ hence as a form of the uncertainty
principle. The conjecture can also be thought of as saying that the sharp
pointwise bounds obtained from the reproducing kernel are the same as one
would expect from the isometry theorem. In the conjecture, $\rho_{t,x}$ is the
heat kernel at a point $x$ in $U/K,$ analytically continued to $U_{\mathbb{C}%
}/K_{\mathbb{C}},$ $\nu_{t}$ is the heat kernel for the dual noncompact
symmetric space $G/K,$ and $j$ is the Jacobian of the exponential mapping for
$G/K.$

\begin{conjecture}
[Heat Kernel Duality]For any simply connected compact symmetric space, there
exist constants $a_{t}$ and $b_{t},$ tending to 1 as $t$ tends to 0, such that%
\[
a_{t}(2\pi t)^{-d}\leq\rho_{t,x}(\exp_{x}(iY))\nu_{t}(Y)j(Y)\leq b_{t}(2\pi
t)^{-d}.
\]

\end{conjecture}

This result is established in \cite{H3} in the compact group case. For
applications to the isometry theorem, one would use this conjecture with $t$
replaced by $2t$ and $Y$ by $2Y.$ (Compare Theorem \ref{sdbnd.thm}.) The paper
\cite{HSt} establishes \textquotedblleft half\textquotedblright\ of this
conjecture (the appropriate bounds on $\nu_{t}$) in the rank-one case and the
even multiplicity case.

\section{Noncompact symmetric spaces\label{noncompact.sec}}

Symmetric spaces (by which I mean, more specifically, Riemannian globally
symmetric spaces) come in three basic types: the compact type, the Euclidean
type, and the noncompact type. Every simply connected symmetric space can be
decomposed as a product of one space of each of these three types. In the
simply connected case, the Euclidean type corresponds simply to manifolds
isometric to $\mathbb{R}^{d}$ with the standard metric. In the simply
connected case, \textquotedblleft compact type\textquotedblright\ is
equivalent to compact, and this is the case we dealt with in Section
\ref{compact.sec}. This leaves, then, the noncompact type. This case is much
trickier, and is just starting to be understood \cite{HM3,HM4,KS1,KS2,KOS}.

In this section, we content ourselves with considering the simplest case,
namely, hyperbolic 3-space $H^{3}$. This is just the noncompact dual (in the
sense described in Section \ref{compact.sec}) of the $S^{3}$ case. The
3-dimensional space is special because this is the only dimension in which
hyperbolic $d$-space is a symmetric space of the \textquotedblleft complex
type,\textquotedblright\ where the complex-type spaces are precisely the
noncompact duals of compact Lie groups. Another way of putting it is that 3 is
the dimension for which the unique positive root for hyperbolic $d$-space has
multiplicity 2. The results described in this section apply more generally to
any noncompact symmetric space of the complex type \cite{HM3,HM4}. I begin by
describing the results and then consider the issue of what the correct
metatheorem might be.

Now, we have already made use in the compact case of the duality between
compact and noncompact symmetric spaces. For spheres, duality is reflected in
the appearance of the heat kernel for hyperbolic $d$-space (the dual to
$S^{d}$) in the isometry and inversion formulas. Roughly, we have spherical
geometry in the real directions of $S_{\mathbb{C}}^{d}$ and hyperbolic
geometry in the imaginary directions. The proof of the inversion formula (and
thus, indirectly, of the isometry formula) hinges on the interplay between
these two geometries.

Now, duality is a symmetric relationship; that is, the dual of the dual is the
original symmetric space. (In the notation of Section \ref{compact.sec}, this
is because multiplying $\mathfrak{p}$ by $i$ twice brings us back to
$\mathfrak{p}$ again.) Thus (in the simplest case), it is reasonable to
attempt to reverse the roles of $S^{3}$ and $H^{3}.$ This means that we now
want to let the base manifold be $H^{3}$, identify the fibers in the tangent
bundle with $S^{3},$ and aim for isometry and inversion formulas using the
heat kernel for $S^{3}$ in the fibers. Unfortunately, further consideration
reveals serious problems with this idea. First, the fibers in $T(H^{3})$ are
diffeomorphic to $\mathbb{R}^{3}$ and not to $S^{3}.$ Second, the well-known
formula for the heat kernel on $H^{3}$ (see (\ref{h3.rhot}) below) reveals
that the analytic continuation of this heat kernel develops singularities once
one travels a sufficient distance in the imaginary direction. The functions
$F$ in the range of the heat operator are going to inherit these singularities
and thus the expression $F(\exp_{x}(iY))$ is going to be defined only for
small $Y.$

These problems can be resolved by making use of certain delicate cancellations
of singularities. Before stating the results, let us set up some notation and
see more precisely what problems arise. We express $H^{3}$ as $G/K$ where $G$
is the identity component of $SO(3,1)$ and where $K$ is $SO(3).$ Then we
define the complexification $H_{\mathbb{C}}^{3}$ of $H^{3}$ to be the complex
manifold $G_{\mathbb{C}}/K_{\mathbb{C}},$ where $G_{\mathbb{C}}$ is
$SO(3,1;\mathbb{C})$ and where $K_{\mathbb{C}}$ is $SO(3;\mathbb{C}).$ Since
$SO(3,1;\mathbb{C})$ is isomorphic to $SO(4;\mathbb{C}),$ $H_{\mathbb{C}}^{3}$
is biholomorphic to $S_{\mathbb{C}}^{3}=SO(4;\mathbb{C})/SO(3;\mathbb{C}).$

We define a map $\Phi$ from $T(H^{3})$ into $H_{\mathbb{C}}^{3}$ by the same
formula as in the sphere case,%
\[
\Phi(x,Y)=\exp_{x}(iY),
\]
where $\exp_{x}(iY)$ refers to the analytic continuation of the geometric
exponential map for $H^{3}.$ The map $\Phi$ is a well-defined smooth map of
$T(H^{3})$ into $H_{\mathbb{C}}^{3},$ but it is not a global diffeomorphism,
nor even a local diffeomorphism. The differential of $\Phi$ is nonsingular
when $\left\vert Y\right\vert <\pi/2,$ but is singular when $\left\vert
Y\right\vert =\pi/2.$ The map $\Phi$ is a diffeomorphism of the set%
\[
T^{\pi/2}(H^{3}):=\left\{  \left.  (x,Y)\in T(H^{3})\right\vert ~\left\vert
Y\right\vert <\pi/2\right\}
\]
with its image in $H_{\mathbb{C}}^{3}$; this image is the so-called
\textquotedblleft crown domain\textquotedblright\ or Akhiezer--Gindikin domain
inside $H_{\mathbb{C}}^{3}.$ If one pulls the complex structure on
$H_{\mathbb{C}}^{3}$ back to $T^{\pi/2}(H^{3})$ by means of $\Phi,$ the result
is the \textquotedblleft adapted complex structure\textquotedblright\ on
$T^{\pi/2}(H^{3}),$ as defined by Guillemin--Stenzel and Lempert--Sz\H{o}ke
\cite{GS1,GS2,LS,Sz}. The adapted complex structure on $T^{\pi/2}(H^{3})$
degenerates as one approaches the boundary, as one would expect from the
degeneration of the map $\Phi.$

Recall that we wish to identify each fiber in $T(H^{3})$ with $S^{3}.$
Although this does not make sense globally, it does make sense locally, using
the exponential map for $S^{3}.$ Under the exponential map, the ball of radius
$\pi/2$ in each fiber is identified with one hemisphere of $S^{3}. $

Kr\"{o}tz and Stanton show \cite[Thm. 6.1]{KS2} that for any $f$ in
$L^{2}(H^{3},dx)$ (where $dx$ is the Riemannian volume measure),
$e^{t\Delta/2}f$ has an analytic continuation to the crown domain in
$H_{\mathbb{C}}^{3}.$ This means that the expression $F(\exp_{x}(iY))$ is
well-defined whenever $\left\vert Y\right\vert <\pi/2$ (but will generally
develop singularities if one tries to go beyond $\left\vert Y\right\vert
=\pi/2$). It is not hard to show, using the explicit form of the spherical
functions on $H^{3},$ a better result for radial functions. Suppose $f\in
L^{2}(H^{3})$ is radial with respect to some basepoint $x_{0}.$ (Since $H^{3}$
is a rank-one symmetric space, this means simply that $f$ is a function of the
distance from $x_{0}.$) Then the map $Y\rightarrow F(\exp_{x_{0}}(iY))$ has a
real analytic extension from the ball of radius $\pi/2$ to the ball of radius
$\pi.$ (Compare Section 4 of \cite{KS1}.) That is, in the radial case,
\textit{when working at the basepoint}, one does not encounter singularities
until $\left\vert Y\right\vert =\pi.$

Finally, we write down the formula \cite[Prop. 3.2]{Ga} for the heat kernel on
$H^{3}$:
\begin{equation}
\nu_{x,t}(\exp_{x}Y)=e^{-t/2}\frac{e^{-\left\vert Y\right\vert ^{2}/2t}}{(2\pi
t)^{3/2}}\frac{\left\vert Y\right\vert }{\sinh\left\vert Y\right\vert
}.\label{h3.rhot}%
\end{equation}
(This heat kernel appeared previously in the isometry and inversion formula
for $S^{3},$ but is now playing a different role.) The factor of $\left\vert
Y\right\vert /\sinh\left\vert Y\right\vert $ is $j^{-1/2}(Y),$ where $j(Y)$ is
the Jacobian of the exponential mapping $\exp_{x}$ at the point $Y.$ If we
evaluate $\nu_{x,t}$ at $\exp_{x}(iY),$ the factor of $\sinh\left\vert
Y\right\vert $ becomes $\sin\left\vert Y\right\vert $ and we encounter the
first singularity at $\left\vert Y\right\vert =\pi.$ This reflects the
expected behavior of radial functions in the range of the heat operator, as
described in the previous paragraph. If, however, one works at some point
other than the basepoint, one will encounter singularities sooner. That is, if
one looks at $\nu_{x_{1},t}(\exp_{x_{2}}(iY))$ with $x_{1}\neq x_{2},$ one
will encounter singularities between $\left\vert Y\right\vert =\pi/2$ and
$\left\vert Y\right\vert =\pi.$ As $x_{1}$ and $x_{2}$ vary, these
singularities will come arbitrarily close to $\left\vert Y\right\vert =\pi/2.$
Thus the analytic continuation of the heat kernel makes sense only on the tube
$T^{\pi/2}(H^{3})$ and not (as it might superficially appear) on the tube
$T^{\pi}(H^{3}).$

We are now ready to think about how one might characterize the range of the
time-$t$ heat operator for $H^{3}.$ In light of the Kr\"{o}tz--Stanton result,
the image of $L^{2}(H^{3},dx)$ under $e^{t\Delta/2}$ is \textit{some} space of
holomorphic functions on $T^{\pi/2}(H^{3}).$ Unfortunately, there does not
seem to be any condition, involving the behavior of $F(\exp_{x}(iY))$ on
$T^{\pi/2}(H^{3}),$ that could reasonably characterize this space. After all,
if we start with $f$ in $L^{2}(H^{3},dx),$ we already get an extension to
$T^{\pi/2}(H^{3})$ as soon as we apply the time-$\varepsilon$ heat operator
for any small, positive $\varepsilon.$ Applying the heat operator for some
larger time $t$ does not give an extension to a larger domain. Furthermore,
the behavior of $F(\exp_{x}(iY))$ as we approach $\left\vert Y\right\vert
=\pi/2$ does not seem to depend on $t.$ Neither the location nor the type of
singularity we get in the analytically continued heat kernel depends on $t.$
On the other hand, the range of the time-$t$ heat operator depends very
strongly on $t.$ It is simply not evident what sort of $t$-dependent condition
we could impose on the values of $F(\exp_{x}(iY))$ for $\left\vert
Y\right\vert <\pi/2$ that could conceivably characterize the image of the
time-$t$ heat operator. (Compare Remark 3.1 in \cite{KOS}.)

It seems, then, that to make progress we must go beyond the set $T^{\pi
/2}(H^{3})$ and find a way to work on the whole tangent bundle. Of course, as
soon as we leave the safe waters of $T^{\pi/2}(H^{3}),$ we encounter
singularities, both in the complex structure and in the functions of the form
$F(\exp_{x}(iY)).$ Nevertheless, I believe that to get a reasonable
characterization of the range of the heat operator, we must say,
\textquotedblleft Damn the singularities! Full speed ahead!\textquotedblright%
\footnote{I am inspired here by the American naval officer, David Glasgow
Farragut, whose strategy of \textquotedblleft Damn the torpedoes! Full speed
ahead!\textquotedblright\ led to victory over Confederate forces in the Battle
of Mobile Bay in 1864.} That is, we must move into the dangerous waters of
$\left\vert Y\right\vert \geq\pi/2$ and simply find some way to deal with the singularities.

The key to dealing with the singularities is to recognize that they have a
\textquotedblleft universal\textquotedblright\ character. This means that
there is a singular quantity that can be factored out such that what remains
is nonsingular and such that the singular factor is \textquotedblleft
universal\textquotedblright---independent of $F$ and $t.$ This is simplest to
explain in the radial case; the general case is more subtle, but in the same
spirit. Note, to start, that the singularity in the analytically continued
heat kernel for $H^{3}$ is independent of $t$; for all $t,$ we have the same
factor of $\left\vert Y\right\vert /\sinh\left\vert Y\right\vert $ that
becomes $\left\vert Y\right\vert /\sin\left\vert Y\right\vert $ upon analytic
continuation. The same behavior is exhibited by general radial functions in
the range of the heat operator. Specifically, if $f\in L^{2}(H^{3},dx)$ is
radial with respect to a basepoint $x_{0}$ and $F=e^{t\Delta/2}f,$ then
\cite[Sect. 3]{HM3} the map%
\[
Y\rightarrow F(\exp_{x_{0}}(Y))\frac{\sinh\left\vert Y\right\vert }{\left\vert
Y\right\vert },\quad Y\in T_{x_{0}}(H^{3}),
\]
has an entire holomorphic extension to the complexified tangent space at
$x_{0}.$ That is to say, if we factor out from $F$ the universal factor of
$\left\vert Y\right\vert /\sinh\left\vert Y\right\vert ,$ what remains (given
above) has no singularities when expressed in exponential coordinates at the
basepoint. This observation leads \cite[Sect. 3]{HM3} to a simple isometry
theorem characterizing the image under $e^{t\Delta/2}$ of the radial functions
in $L^{2}(H^{3},dx).$

In the nonradial case, the singularities in $F$ itself are, unfortunately, not
of a universal character. Nevertheless, the inversion and isometry theorems
can be expressed in terms of certain integrals involving $F,$ and these
integral turn out to have universal singularities. Without further ado, let us
state the inversion and isometry results for $H^{3}.$ The inversion formula is
proved in \cite[Sect. 4]{HM3}. The isometry formula will be addressed in
\cite{HM4}; since the details have not all been written down yet, I\ state the
isometry formula as a conjecture. Another, seemingly different isometry
formula (for general symmetric spaces of the noncompact type) has been
obtained independently by Kr\"{o}tz, \'{O}lafsson, and Stanton \cite{KOS}. It
remains to be worked out how this formula relates to the one in Conjecture
\ref{noncompact.conj}

\begin{theorem}
\label{h3inv.thm}For any $f\in L^{2}(H^{3}),$ the function%
\begin{equation}
L(x,R):=\int_{\substack{Y\in T_{x}(H^{3}) \\\left\vert Y\right\vert \leq R
}}F(\exp_{x}(iY))e^{t/2}\frac{e^{-\left\vert Y\right\vert ^{2}/2t}}{(2\pi
t)^{3/2}}\frac{\sin\left\vert Y\right\vert }{\left\vert Y\right\vert
}~dY,\label{h3inv.inta}%
\end{equation}
initially defined for $R$ in the interval $(0,\pi/2),$ has a real-analytic
extension to the interval $(0,\infty).$ If $f$ is sufficiently smooth, then
\[
f(x)=\lim_{R\rightarrow\infty}L(x,R)
\]
and we may write, informally,%
\begin{equation}
f(x)=\text{\textquotedblleft}\lim_{R\rightarrow\infty}\text{{}%
\textquotedblright}\int_{\substack{Y\in T_{x}(H^{3}) \\\left\vert Y\right\vert
\leq R}}F(\exp_{x}(iY))e^{t/2}\frac{e^{-\left\vert Y\right\vert ^{2}/2t}%
}{(2\pi t)^{3/2}}\frac{\sin\left\vert Y\right\vert }{\left\vert Y\right\vert
}~dY.\label{h3inv.int}%
\end{equation}

\end{theorem}

\begin{conjecture}
\label{noncompact.conj}For any $f\in L^{2}(H^{3}),$ the quantity%
\begin{equation}
M(R):=\int_{x\in H^{3}}\int_{\substack{Y\in T_{x}(H^{3})\\\left\vert
Y\right\vert \leq R}}\left\vert F\left(  \exp_{x}\left(  iY\right)  \right)
\right\vert ^{2}e^{t}\frac{\sin\left\vert 2Y\right\vert }{\left\vert
2Y\right\vert }\frac{e^{-\left\vert Y\right\vert ^{2}/t}}{(\pi t)^{3/2}%
}~dY~dx\label{mr}%
\end{equation}
is finite for all $R<\pi/2,$ where $dY$ denotes Lebesgue measure on
$T_{x}(H^{3})$ and where $dx$ denotes the Riemannian volume measure on
$H^{3}.$ Furthermore, $M$ has a real-analytic continuation from $(0,\pi/2)$ to
$(0,\infty)$ that satisfies%
\[
\lim_{R\rightarrow\infty}M(R)=\int_{H^{3}}\left\vert f(x)\right\vert ^{2}~dx.
\]
Thus we may write, informally,%
\begin{align}
& \int_{H^{3}}\left\vert f(x)\right\vert ^{2}~dx\label{h3isom.int}\\
& =\text{\textquotedblleft}\lim_{R\rightarrow\infty}\text{{}\textquotedblright%
}\int_{x\in H^{3}}\int_{\substack{Y\in T_{x}(H^{3})\\\left\vert Y\right\vert
\leq R}}\left\vert F\left(  \exp_{x}\left(  iY\right)  \right)  \right\vert
^{2}e^{t}\frac{\sin\left\vert 2Y\right\vert }{\left\vert 2Y\right\vert }%
\frac{e^{-\left\vert Y\right\vert ^{2}/t}}{(\pi t)^{3/2}}~dY~dx.\nonumber
\end{align}

Conversely, suppose that $F$ is a function on $H^{3}$ having an analytic
continuation to $T^{\pi/2}(H^{3})$, with respect to the adapted complex
structure. Suppose further that the integral $M(R)$ defined in (\ref{mr}) is
finite for all $R<\pi/2,$ that $M$ has an analytic continuation to
$(0,\infty),$ and that $\lim_{R\rightarrow\infty}M(R)$ exists and is finite.
Then there exists a unique $f$ in $L^{2}(H^{3},dx)$ for which $F=e^{t\Delta
/2}f.$
\end{conjecture}

Observe how closely these formulas parallel the corresponding results for
$S^{3}$ in Theorem \ref{s3.thm}. To get from (\ref{h3inv.int}) and
(\ref{h3isom.int}) to the corresponding $S^{3}$ versions (Theorem
\ref{s3.thm}), we simply change $\sin$ to $\sinh$; change the exponential
factors $e^{t/2} $ and $e^{t}$ to $e^{-t/2}$ and $e^{-t},$ respectively; and
omit the quotation marks in the limit as $R$ tends to infinity. Furthermore,
the quantity%
\[
e^{t/2}\frac{e^{-\left\vert Y\right\vert ^{2}/2t}}{(2\pi t)^{3/2}}\frac
{\sin\left\vert Y\right\vert }{\left\vert Y\right\vert }~dY
\]
in (\ref{h3inv.int}) and (scaled appropriately) in (\ref{h3isom.int}) is
essentially just the heat kernel measure for $S^{3}.$ More precisely, this is
an \textquotedblleft unwrapped\textquotedblright\ version of the heat kernel
measure for $S^{3},$ in a sense described in \cite{HM3}.

Thus, if we are sufficiently imaginative, we can achieve something close to
the suggestion made earlier in this section: We let the base manifold be
$H^{3},$ identify the fibers with $S^{3},$ and obtain isometry and inversion
formulas involving the heat kernel for $S^{3}$ in the fibers. (It is slightly
more accurate to say that we identify the fibers with the \textit{tangent
space} to $S^{3}$ at the basepoint, where this tangent space is then
identified \textit{locally} with $S^{3}$ itself.)

The idea behind the analytic continuation of the function $L(x,R)$ in Theorem
\ref{h3inv.thm} is as follows. The integral in (\ref{h3inv.inta}) only
\textquotedblleft sees\textquotedblright\ the radial part of the function
$Y\rightarrow F(\exp_{x}(iY)).$ This radial part, as discussed earlier, has
singularities only of the $\left\vert Y\right\vert /\sin\left\vert
Y\right\vert $ variety, and these are canceled by the factor of $\sin
\left\vert Y\right\vert /\left\vert Y\right\vert $ in (\ref{h3inv.inta}).

Let us see if we can formulate a metatheorem that would capture the spirit of
these two results. Note that in the sphere case, we actually have two
different versions of the metatheorem. There is the one in Metatheorem
\ref{sphere.meta} and the one in (\ref{approx.bnd}). In the $d=3$ case, the
first of these conditions will involve $\left\vert Y\right\vert /\sinh
\left\vert Y\right\vert $ and the second $(\left\vert 2Y\right\vert
/\sinh\left\vert 2Y\right\vert )^{1/2}.$ The first condition is the one
suggested by the inversion formula and by the growth (\ref{rhot.bnd}) of the
analytically continued heat kernel. The second condition is the one suggested
by the isometry formula and the sharp pointwise bounds (Theorem
\ref{sdbnd.thm} or (\ref{precise.bnd})). In the sphere case, these two
conditions differ only by a polynomial factor and we do not really need to
distinguish them.

In the $H^{3}$ case, we could consider the naive hyperbolic analogs of these
two conditions, obtained by changing hyperbolic sine to ordinary sine. This
means the condition on growth of $F(\exp_{x}(iY))$ should be either growth at
most like%
\begin{equation}
e^{\left\vert Y\right\vert ^{2}/2t}\frac{\left\vert Y\right\vert }%
{\sin\left\vert Y\right\vert }\label{inv.sing}%
\end{equation}
or at most like%
\begin{equation}
e^{\left\vert Y\right\vert ^{2}/2t}\left(  \frac{\left\vert 2Y\right\vert
}{\sin\left\vert 2Y\right\vert }\right)  ^{1/2}.\label{inv.isom}%
\end{equation}

In the $H^{3}$ case, these conditions are seemingly quite different, since
they suggest different things about the locations of the singularities.
Nevertheless, both conditions suggest reasonable metatheorems; it is just that
the two conditions correspond to two different ways of canceling out the
singularities. This leads us, then, to two different (hopefully approximately
equivalent!) candidates for a metatheorem in the $H^{3}$ case.

\begin{metatheorem}
[Version 1]\label{hyper.meta1}The functions in the range of the time-$t$ heat
operator $e^{t\Delta/2}$ for $H^{3}$ are those functions $F$ such that for
each $x,$ the map%
\[
Y\rightarrow F(\exp_{x}(iY))^{\mathrm{rad}}\frac{\sin\left\vert Y\right\vert
}{\left\vert Y\right\vert }%
\]
has a real-analytic extension from $\left\{  \left\vert Y\right\vert
<\pi/2\right\}  $ to all of $T_{x}(H^{3})$ with growth at most like
$e^{\left\vert Y\right\vert ^{2}/2t}.$
\end{metatheorem}

Here $F(\exp_{x}(iY))^{\mathrm{rad}}$ means the radial part of the function
$F(\exp_{x}(iY)),$ with respect to $Y$ with $x$ fixed, that is, the average of
$F(\exp_{x}(iY))$ with respect to the action of the rotation group on the $Y$ variable.

\begin{metatheorem}
[Version 2]\label{hyper.meta2}The functions in the range of the time-$t$ heat
operator $e^{t\Delta/2}$ for $H^{3}$ are those functions $F$ such that the map%
\begin{equation}
R\rightarrow\frac{\sin(2R)}{2R}\int_{x\in H^{3}}\int_{\substack{Y\in
T_{x}(H^{3}) \\\left\vert Y\right\vert =R}}\left\vert F(\exp_{x}%
(iY)\right\vert ^{2}~dY~\label{meta.int2}%
\end{equation}
is finite for all $R<\pi/2$ and has a real-analytic continuation from
$R\in(0,\pi/2)$ to \thinspace$(0,\infty)$ with growth at most like
$e^{R^{2}/t}.$
\end{metatheorem}

Metatheorem \ref{hyper.meta2} should be interpreted carefully. After all, if
$f$ is a nice function that just happens not to be in $L^{2}(H^{3},dx),$ then
$F$ may not be in $L^{2}(H^{3},dx)$ either, and in that case the integral in
(\ref{meta.int2}) is likely to be infinite even for small $R.$ Thus
Metatheorem \ref{hyper.meta2} is likely to be most useful for classes of
initial conditions where $e^{\varepsilon\Delta/2}f$ is in $L^{2}(H^{3},dx)$
for all positive $\varepsilon.$ Such classes include positive and negative
Sobolev spaces and $L^{p}$ spaces for $1\leq p\leq2.$ For example, Metatheorem
\ref{hyper.meta2} suggests a possible characterization of the image under
$e^{t\Delta/2}$ of the intersection over $n$ of the $n^{\text{th}}$ Sobolev
space on $H^{3}.$ On the other hand, Metatheorem \ref{hyper.meta1} is more
flexible an ought to apply to classes of initial conditions with varying rates
of growth or decay at infinity.

I conclude this section with a discussion of multiplication properties in the
$H^{3}$ case. The sorts of results we had (Theorems \ref{mult.thm1},
\ref{mult.thm2}, and \ref{mult.thm3}) in the case of Euclidean space or a
compact symmetric space do \textit{not} extend to the case of hyperbolic
3-space. Of course, if $F_{1}$ and $F_{2}$ are both holomorphic on $T^{\pi
/2}(H^{3})$ (with respect to the complex structure obtained by using $\Phi$),
then so is $F_{1}F_{2}.$ However, the cancellation of singularities we have
outside of $T^{\pi/2}(H^{3})$ will not be preserved under multiplication. This
is seen most easily in the radial case. If $F_{1}$ and $F_{2}$ are radial
functions in the image of the heat operator, then $F_{1}(\exp_{x_{0}}(iY))$
and $F_{2}(\exp_{x_{0}}(iY))$ will both have singularities like $1/\sin
\left\vert Y\right\vert $ as we approach $\left\vert Y\right\vert =\pi.$ When
we multiply them, we will get a singularity like $1/\sin^{2}\left\vert
Y\right\vert $. A function with a singularity of this sort cannot (in light of
Section 3 of \cite{HM3}) be in the range of the heat operator for any positive time.

In the nonradial case, we can see that neither Metatheorem \ref{hyper.meta1}
nor Metatheorem \ref{hyper.meta2} describes a class of functions that will be
closed under multiplication, even if we permit a change in the time parameter.
In Metatheorem \ref{hyper.meta1}, the first problem is that the radialization
of a product is not the product of the radializations, and the second problem
is that even when the functions are already radial, the factor of
$\sin|Y|/|Y|$ can only cancel the singularities in one of the two functions.
In Metatheorem \ref{hyper.meta2}, we expect that the divergence of the
integral as $R$ approaches $\pi/2$ to be worse for $F_{1}F_{2}$ than for
either function alone. Thus, the divergence of the integral for $F_{1}F_{2}$
is probably not going to be canceled by the factor of $\sin(2R)/2R.$

The results of Kr\"{o}tz and Stanton in \cite{KS1} suggest that $F_{1}F_{2}$
is regular enough to be in the range of the operator $\exp[-(\pi
/2)\sqrt{-\Delta}],$ but no better than this.

\section{Concluding remarks}

I conclude this paper by discussing what little is known (from the point of
view of this paper) about the range of the heat operator for manifolds other
than symmetric spaces. A recent paper of Kr\"{o}tz, Thangavelu, and Xu
\cite{KTX} calculates the image of the heat operator for a left-invariant
Riemannian (\textit{not} sub-Riemannian) metric on the Heisenberg group. They
characterize this image as a certain space of holomorphic functions on the
complexified Heisenberg group. However, the description of this space is a bit
complicated; it is the sum of two weighted Bergman spaces with weight
functions that assume both positive and negative values.

For a general compact real analytic manifold $\mathcal{M},$ there is a
characterization of the range of the Poisson semigroup $\exp[-t\sqrt{-\Delta
}]$ that gives an idea of how tricky the range of the heat operator is likely
to be in general. According to \cite{GS1,GS2,LS,Sz}, there exists a
canonically defined \textquotedblleft adapted complex
structure\textquotedblright\ on a sufficiently small tube $T^{\varepsilon
}(\mathcal{M})$ inside the tangent bundle. Then a result of \cite{GS2}, based
in part on earlier work of L. Boutet de Monvel, is as follows.

\begin{theorem}
There exists some $\varepsilon_{0}\leq\varepsilon$ such that for all
$t<\varepsilon_{0},$ a function $F$ on $\mathcal{M}$ is of the form
$F=\exp[-t\sqrt{-\Delta}]f,$ with $f$ in $C^{\infty}(\mathcal{M}),$ if and
only if $F$ has an analytic continuation to the tube $T^{t}(\mathcal{M})$ that
is smooth up to the boundary.
\end{theorem}

That is, roughly, the range of the operator $\exp[-t\sqrt{-\Delta}]$ is the
space of functions having an analytic continuation to a tube of radius $t,$
provided that $t<\varepsilon_{0}.$Note that this theorem says nothing about
the range of $\exp[-t\sqrt{-\Delta}]$ when $t>\varepsilon_{0}.$ Also, it is
not known whether it is possible to take $\varepsilon_{0}$ to equal
$\varepsilon$ (the largest radius for which the adapted complex structure
exists), as would be natural to expect.

This result indicates the challenges that await in trying to characterize the
range of the heat operator. If $F$ is in the range of $e^{\varepsilon\Delta}$
for \textit{any} positive $\varepsilon,$ then $F$ is automatically in the
range of $\exp[-t\sqrt{-\Delta}]$ for \textit{all} positive $t.$ This means
that even for very small $t,$ a function $F$ range of the time-$t$ heat
operator will be holomorphic on the tube of radius $\varepsilon_{0}.$ Quite
possibly, $F$ is holomorphic on the tube of radius $\varepsilon,$ the largest
tube on which the adapted complex structure exists. But this condition is
unlikely to be sufficient, since holomorphicity on a tube of finite radius
characterizes the range of the Poisson semigroup and not the heat semigroup.
Conceivably, the way $F$ behaves as on approaches the boundary of the tube of
radius $\varepsilon$ could characterize the range of the heat operator, but
our experience in the $H^{3}$ case make this seem unlikely. More likely, if
there is a characterization of the range of the heat operator, it will have to
involve the nature of the singularities of $F$ once one moves beyond the
maximal tube on which the adapted complex structure is defined.

Manifolds where the adapted complex structure exists on all of $T(\mathcal{M}%
)$ might be more tractable to deal with, but these appear to be rare. (Compare
\cite{A}.)


\begin{thebibliography}{99999}                                                                                            %
\bibitem[A]{A}R. M. Aguilar, Symplectic reduction and the homogeneous complex
Monge-Amp\`{e}re equation, \textit{Ann. Global Anal. Geom.} \textbf{19}
(2001), 327--353.

\bibitem[BSZ]{BSZ}J. C. Baez, I. E. Segal, and Z.-F. Zhou, Introduction to
algebraic and constructive quantum field theory. Princeton Series in Physics.
Princeton University Press, Princeton, NJ, 1992.

\bibitem[B1]{Ba1}V. Bargmann, On a Hilbert space of analytic functions and an
associated integral transform, \textit{Comm. Pure Appl. Math.} \textbf{14}
(1961), 187--214.

\bibitem[B2]{Ba2}V. Bargmann, On a Hilbert space of analytic functions and an
associated integral transform. Part II. A family of related function spaces.
Application to distribution theory, \textit{Comm. Pure Appl. Math.}
\textbf{20} (1967), 1--101.

\bibitem[Da]{Da}E. B. Davies, Heat kernels and spectral theory. Cambridge
Tracts in Mathematics, 92. Cambridge University Press, Cambridge, 1990.

\bibitem[DM]{DM}E. B. Davies and N. Mandouvalos, Heat kernel bounds on
hyperbolic space and Kleinian groups, \textit{Proc. London Math. Soc.} (3)
\textbf{57} (1988), 182--208.

\bibitem[D]{D}B. K. Driver, On the Kakutani-It\^{o}-Segal-Gross and
Segal-Bargmann-Hall isomorphisms, \textit{J. Funct. Anal.} \textbf{133}
(1995), 69--128.

\bibitem[DG]{DG}B. K. Driver and L. Gross, Hilbert spaces of holomorphic
functions on complex Lie groups. In: New trends in stochastic analysis (K.
Elworthy, S. Kusuoka, and I. Shigekawa, Eds.), 76--106, World Sci. Publishing,
River Edge, NJ, 1997.

\bibitem[DH]{DH1}B. K. Driver and B. C. Hall, Yang--Mills theory and the
Segal--Bargmann transform, \textit{Comm. Math. Phys.} \textbf{201} (1999), 249--290.

\bibitem[FMMN1]{FMMN}C. Florentino, P. Matias, J. Mour\~{a}o, and J. Nunes,
Geometric quantization, complex structures, and the coherent state transform,
\textit{J. Funct. Anal.} \textbf{221} (2005), 303--322.

\bibitem[FMMN2]{FMMN2}C. Florentino, P. Matias, J. Mour\~{a}o, and J. Nunes,
On the BKS pairing for Kahler quantizations of the cotangent bundle of a Lie
group, preprint.

http://arxiv.org/abs/math.DG/0411334

\bibitem[Ga]{Ga}R. Gangolli, Asymptotic behavior of spectra of compact
quotients of certain symmetric spaces, \textit{Acta Math.} \textbf{121}
(1968), 151--192

\bibitem[Gr]{Gr}L. Gross, Uniqueness of ground states for Schr\"{o}dinger
operators over loop groups, \textit{J. Funct. Anal.} \textbf{112} (1993), 373--441.

\bibitem[GS1]{GS1}V. Guillemin and M. B. Stenzel, Grauert tubes and the
homogeneous Monge-Amp\`{e}re equation, \textit{J. Differential Geom.}
\textbf{34} (1991), 561--570.

\bibitem[GS2]{GS2}V. Guillemin and M. B. Stenzel, Grauert tubes and the
homogeneous Monge-Amp\`{e}re equation. II, \textit{J. Differential Geom.}
\textbf{35} (1992), 627--641.

\bibitem[H1]{H1}B. C. Hall, The Segal-Bargmann ''coherent state'' transform
for compact Lie groups, \textit{J. Funct. Anal.} \textbf{122} (1994), 103--151.

\bibitem[H2]{H2}B. C. Hall, The inverse Segal-Bargmann transform for compact
Lie groups, \textit{J. Funct. Anal.} \textbf{143} (1997), 98--116.

\bibitem[H3]{H3}B. C. Hall, Phase space bounds for quantum mechanics on a
compact Lie group, \textit{Comm. Math. Phys.} \textbf{184} (1997), 233--250.

\bibitem[H5]{newform}B. C. Hall, A new form of the Segal-Bargmann transform
for Lie groups of compact type, \textit{Canad. J. Math.} \textbf{51} (1999), 816--834.

\bibitem[H6]{mexnotes}B. C. Hall, Holomorphic methods in analysis and
mathematical physics. In: First Summer School in Analysis and Mathematical
Physics (S. P erez-Esteva and C. Villegas-Blas, Eds.), 1--59, Contemp. Math.,
260, Amer. Math. Soc., Providence, RI, 2000.

\bibitem[H7]{geoquant}B. C. Hall, Geometric quantization and the generalized
Segal--Bargmann transform for Lie groups of compact type, \textit{Comm. Math.
Phys.} \textbf{226} (2002), 233-268.

\bibitem[H8]{bull}B. C. Hall, Harmonic analysis with respect to heat kernel
measure, \textit{Bull. Amer. Math. Soc. (N.S.)} \textbf{38} (2001), 43--78.

\bibitem[H9]{lpbounds}B. C. Hall, Bounds on the Segal-Bargmann transform of
$L^{p}$ functions, \textit{J. Fourier Anal. Appl.} \textbf{7} (2001), 553--569.

\bibitem[H10]{ergodic}B. C. Hall, The Segal-Bargmann transform and the Gross
ergodicity theorem. In: Finite and infinite dimensional analysis in honor of
Leonard Gross (H.-H. Kuo and A. N. Sengupta, Eds.), 99--116, Contemp. Math.,
317, Amer. Math. Soc., Providence, RI, 2003.

\bibitem[HL]{HL}B. C. Hall and W. Lewkeeratiyutkul, Holomorphic Sobolev spaces
and the generalized Segal--Bargmann transform, \textit{J. Funct. Anal.}
\textbf{217} (2004), 192-220.

\bibitem[HM1]{HM1}B. C. Hall and J. J. Mitchell, Coherent states on spheres,
\textit{J. Math. Phys.} \textbf{43} (2002), 1211--1236.

\bibitem[HM2]{HM2}B. C. Hall and J. J. Mitchell, The large radius limit for
coherent states on spheres. In, \textquotedblleft Mathematical results in
quantum mechanics (Taxco, 2001),\textquotedblright\ 155--162, Contemp. Math.,
307, Amer. Math. Soc., Providence, RI, 2002.

\bibitem[HM3]{HM3}B. C. Hall and J. J. Mitchell, The Segal--Bargmann transform
for noncompact symmetric spaces of the complex type, \textit{J. Funct. Anal.},
to appear.

\bibitem[HM4]{HM4}B. C. Hall and J. J. Mitchell, The Segal--Bargmann transform
for noncompact symmetric spaces II. In preparation.

\bibitem[HS]{HSt}B. C. Hall and M. B. Stenzel, Sharp bounds for the heat
kernel on certain symmetric spaces of non-compact type. In: Finite and
infinite dimensional analysis in honor of Leonard Gross (H.-H. Kuo and A. N.
Sengupta), 117--135, Contemp. Math., 317, Amer. Math. Soc., Providence, RI, 2003.

\bibitem[He]{He}S. Helgason, \textquotedblleft Differential Geometry, Lie
Groups, and Symmetric Spaces.\textquotedblright\ Academic Press, 1978.

\bibitem[KR]{KR2}K. Kowalski and J. Rembieli\'{n}ski, The Bargmann
representation for the quantum mechanics on a sphere, \textit{J. Math. Phys.}
\textbf{42} (2001), 4138--4147.

\bibitem[KOS]{KOS}B. Kr\"{o}tz, G. \'{O}lafsson, and R. Stanton, The image of
the heat kernel transform on Riemannian symmetric spaces of the non-compact
type, preprint.

http://arxiv.org/abs/math.CA/0407391

\bibitem[KS1]{KS1}B. Kr\"{o}tz and R. J. Stanton, Holomorpic extensions of
representations: (I) automorphic functions, \textit{Ann. Math.} \textbf{159}
(2004), 641-724.

\bibitem[KS2]{KS2}B. Kr\"{o}tz and R. J. Stanton, Holomorphic extension of
representations: (II) geometry and harmonic analysis, preprint.

http://arxiv.org/abs/math.RT/0502411

\bibitem[KTX]{KTX}B. Kr\"{o}tz, S. Thangavelu, and Y. Xu, The heat kernel
transform for the Heisenberg group, \textit{J. Funct. Anal.}, to appear.

\bibitem[LGS]{LGS}E. Leichtnam, F. Golse, and M. B. Stenzel, Intrinsic
microlocal analysis and inversion formulae for the heat equation on compact
real-analytic Riemannian manifolds, \textit{Ann. Sci. \'{E}cole Norm. Sup.}
(4) \textbf{29} (1996), 669--736.

\bibitem[LS]{LS}L. Lempert and R. Sz\H{o}ke, Global solutions of the
homogeneous complex Monge-Amp\`{e}re equation and complex structures on the
tangent bundle of Riemannian manifolds, \textit{Math. Ann.} \textbf{290}
(1991), 689--712.

\bibitem[Se1]{Se1}I. E. Segal, Mathematical problems of relativistic physics.
In: \textquotedblleft Proceedings of the Summer Seminar, Boulder, Colorado,
1960\textquotedblright\ (M. Kac, Ed.), American Mathematical Society,
Providence, RI, 1963.

\bibitem[Se2]{Se2}I. E. Segal, Mathematical characterization of the physical
vacuum for a linear Bose-Einstein field. (Foundations of the dynamics of
infinite systems. III) \textit{Illinois J. Math.} \textbf{6} (1962), 500--523.

\bibitem[Se3]{Se3}I. E. Segal, The complex-wave representation of the free
Boson field. In: \textquotedblleft Topics in Functional
Analysis\textquotedblright\ (I. Gohberg and M. Kac, Eds.), Advances in
Mathematics Supplementary Studies, Vol. 3, Academic Press, New York, 1978.

\bibitem[St]{St}M. B. Stenzel, The Segal-Bargmann transform on a symmetric
space of compact type, \textit{J. Funct. Anal.} \textbf{165} (1999), 44--58.

\bibitem[Sz]{Sz}R. Sz\H{o}ke, Complex structures on tangent bundles of
Riemannian manifolds, \textit{Math. Ann.} \textbf{291} (1991), 409--428.

\bibitem[Ta]{Ta}M. E. Taylor, Noncommutative harmonic analysis. Mathematical
Surveys and Monographs, 22. American Mathematical Society, Providence, RI, 1986.

\bibitem[TW]{Th}T. Thiemann and O. Winkler, Gauge field theory coherent states
(GCS): II. Peakedness properties, \textit{Classical Quantum Gravity}
\textbf{18} (2001), 2561-2636.
\end{thebibliography}
\end{document}